\documentclass[aps,nofootinbib,notitlepage,superscriptaddress]{revtex4-1}
\setcitestyle{authoryear,round,semicolon,aysep={ }}
\usepackage{amsmath}
\usepackage{amssymb}
\usepackage{amsfonts}
\usepackage{graphicx}
\usepackage{hyperref}
\usepackage{color}
\usepackage{multirow}
\usepackage{mathtools}
\usepackage{booktabs}
\usepackage{amsthm}
\usepackage{mathtools}
\usepackage{color}
\usepackage{url}

\hypersetup{pdftitle=Simple Geometry for Networks}

\begin{document}
\title{A Simple Differential Geometry for Complex Networks}

\author{Emil Saucan}
\email{Correspondence to: semil@braude.ac.il}
\affiliation{Department of Applied Mathematics, ORT Braude College, 
Karmiel 2161002, Israel}
\author{Areejit Samal}
\affiliation{The Institute of Mathematical Sciences (IMSc), Homi Bhabha National
Institute (HBNI), Chennai 600113 India}
\affiliation{Max Planck Institute for Mathematics in the Sciences, Leipzig 04103
Germany}
\author{J\"{u}rgen Jost}
\affiliation{Max Planck Institute for Mathematics in the Sciences, Leipzig 04103
Germany}
\affiliation{The Santa Fe Institute, Santa Fe, New Mexico 87501, USA}

\begin{abstract}
We introduce new definitions of sectional, Ricci and scalar curvature for 
networks and their higher dimensional counterparts, derived from two classical 
notions of curvature for curves in general metric spaces, namely, the Menger 
curvature and the Haantjes curvature. These curvatures are applicable to 
unweighted or weighted and undirected or directed networks, and are more 
intuitive and easier to compute than other network curvatures. In particular, 
the proposed curvatures based on the interpretation of Haantjes definition as 
geodesic curvature allow us to give a network analogue of the classical local 
Gauss-Bonnet theorem. Furthermore, we propose even simpler and more intuitive 
proxies for the Haantjes curvature that allow for even faster and easier 
computations in large-scale networks. In addition, we also investigate the 
embedding properties of the proposed Ricci curvatures. Lastly, we also 
investigate the behaviour, both on model and real-world networks, of the 
curvatures introduced herein with more established notions of Ricci curvature 
and other widely-used network measures.
\end{abstract}

\maketitle

\section{Introduction}

Till recently, mathematical analysis of complex networks was largely based 
upon combinatorial invariants and models \cite{Watts1998,Barabasi1999,
Albert2002,Newman2010,Dorogovtsev2013}, at the detriment of the geometric 
approach. Especially, theoretical developments in discrete differential 
geometry leading to notions of curvature \cite{Ollivier2009,Forman2003} for 
graphs or networks remained largely unexplored until recently in network 
science \cite{Ni2015,Sandhu2015,Sreejith2016}. Noteworthy, one of the 
widely-used combinatorial measure, the \textit{clustering coefficient} 
\cite{Watts1998}, to characterize complex networks is in fact, a discretization 
of the classical Gauss curvature \cite{Eckmann2002}. Lately, notions of network 
curvature have proven to be an important tool in the analysis of complex 
networks \cite{Ni2015,Sandhu2015,Sreejith2016,Samal2018}. In particular, 
Ollivier's Ricci curvature \cite{Ollivier2009} has been extensively used in the 
analysis of complex networks in its various avatars \cite{Ni2015,Sandhu2015,
Asoodeh2018,Ni2019,Farooq2019}. Based on Forman's work \cite{Forman2003}, 
another approach towards the introduction of Ricci curvature in the study of 
networks was also proposed \cite{Sreejith2016,Weber2017}. Moreover, the two 
above-mentioned notions of Ricci curvature in several model and real-world 
networks were compared \cite{Samal2018}. However, these two notions of Ricci 
curvature for networks have certain drawbacks. Their common denominator is the 
theoretical advanced apparatus which stands at the base of both notions of 
discrete Ricci curvature for networks. While Ollivier-Ricci curvature is 
prohibitively hard to compute in large networks, Forman-Ricci curvature is 
extremely simple to compute in large networks \cite{Samal2018}. On the other 
hand, Forman-Ricci curvature is less intuitive than Ollivier-Ricci curvature, as 
it is based on a discretization of the so called \textit{Bochner-Weitzenb{\"o}ck 
formula} \cite{Jost2017}. These issues may deter the interdisciplinary community 
of engineers, social scientists and biologists active in network science from 
employing these two notions of discrete Ricci curvature in their research. 

Therefore, it is worthwhile to ask whether it is possible to define more 
intuitive notions of curvature (specifically, Ricci curvature) for networks 
which are modelled using the simple framework of graphs. The answer to this 
quandary is immediate after the realization that the simplest, more general, 
geometric notion at our disposal is that of metric space. Furthermore, 
unweighted and weighted graphs can be easily endowed with a metric, be it the 
combinatorial weights, the path distance, the given weights, the Wasserstein 
metric (also called the Earth mover's distance) \cite{Vaserstein1969}, or a 
metric such as the path degree metric \cite{Dodziuk1986,Keller2015} that can 
be easily obtained from the given weights. For metric spaces, there indeed 
exist a number of simple and intuitive notions of curvature, and in particular, 
two definitions of curvature of curves which go back to Menger \cite{Menger1930} 
and Haantjes \cite{Haantjes1947}, respectively. By partially extending ideas 
previously applied in the context of imaging and graphics ($PL$ manifolds in 
general) \cite{Saucan2015,Gu2013,Saucan2018}, we show that the definitions 
of Menger \cite{Menger1930} and Haantjes \cite{Haantjes1947} allow us to 
naturally define Ricci and scalar curvatures for networks and hypernetworks 
including unweighted, weighted, undirected and directed ones.

The simplest and better known among the two notions of metric curvature is the 
Menger curvature. Notably, the Menger curvature defined for metric triangles 
naturally allows us to define scalar and Ricci curvatures for networks, 
hypernetworks, simplicial complexes and clique complexes including unweighted, 
weighted, undirected and directed ones. Furthermore, we show that this approach 
holds not only for the usual, Euclidean model geometry, but as well for 
hyperbolic and spherical background geometry. We believe this is an important 
feature of the Menger curvature due to the following reason. Currently, while 
network embedding is commonly considered, such embedding is either purely 
combinatorial or, at most, equipped with an undescriptive Euclidean geometry 
which is decoupled from the one of the ambient space, and thus, a vital part 
of the expressiveness of the embedding is lost.

In comparison to the Menger curvature, Haantjes curvature is far less familiar 
among researchers but is a much more flexible notion. This is because Haantjes 
curvature is not restricted solely to triangles, but rather holds for general 
metric arcs. In consequence, Haantjes curvature presents two distinct advantages. 
Firstly, it is applicable to any 2-cell, not just to triangles. Secondly, it can 
be better used as discrete version of the classical geodesic curvature (of curves 
on smooth surfaces). In consequence, Haantjes curvature has a clear advantage 
over Menger curvature as it is applicable to networks or graphs, without 
\textit{any assumption on the background geometry}. In fact, not only can it be
employed in networks of \textit{variable curvature}, it can be used to 
\textit{define} the curvature of such a discrete space. In effect, this allows 
us to study the \textit{intrinsic}, not just the \textit{extrinsic}, geometry of 
networks. From a practical view of applications, the definition of Haantjes 
curvature allows inference of geometric, highly descriptive properties from its 
given characteristic, rather than presume them. From a more theoretical 
perspective but still with practical importance, for instance, while predicting 
the long time behavior, Haantjes curvature enables us to study networks as 
geometric spaces in their own right, not just as graphs realized in some largely 
arbitrary ambient familiar space.

We should emphasize that there is not a single, restricted motivation residing 
behind our quest for a purely metric differential geometry for networks. 
To be sure, the search for faster and more efficient computational tools is 
one of the reasons we adopted this path, but it is not the only one, nor is it, 
by far, the most important one. 

A substantial motivation behind the choice of the metric route towards the 
geometrization of networks stems from the fact that it allows us to consistently 
define and study the {\it intrinsic} geometry of networks, as opposed to the 
more common approach via {\it extrinsic} geometry, that is by embedding the 
networks in suitable ambient space, and studying the geometric properties induced 
by the surrounding space. Lately, this embedding space has usually been taken 
to be the hyperbolic plane or space, given the understanding that, in a sense, 
it combines the efficiency and convenient visualization, mostly using the 
Poincar\'{e} models, see e.g. \cite{Boguna2010,Boguna2020}.

While the later approach is more innate and conducive to intuitive illustrations, 
such embeddings are, apart from the very stringent {\it isometric embeddings}, 
distorting. In contrast, the intrinsic approach, that is the study of networks 
{\it per se}, is independent of any specific embedding, hence of the necessary 
additional computations and also of any distortion. Moreover, this approach allows 
for the independent study of such powerful tools as the Ricci flow, without the 
vagaries associated to the embedding in an ambient space of certain dimension 
(see, e.g. \cite{Saucan2012a} and the bibliography therein). Indeed, the 
development of a purely metric, intrinsic geometry of networks is the driving 
motivation of Section 3.1 and Section 4. 

Our approach to the the development of a metric, intrinsic geometry for networks 
follows the Riemannian paradigm to geometry, more specifically its $PL$ (and 
polyhedral) avatar. In particular, we follow the classical approach to $PL$ 
Geometry (and Topology), in particular the pioneering works of Stone 
\cite{Stone1973,Stone1976}, on the discretization of sectional and Ricci curvatures. 
See also \cite{Gu2013} for a more recent metric extension of Stone's ideas. 
In this approach, polygons (i.e. the 2-dimensional elements in a polyhedral manifold) 
incident to an edge are the direct and natural $PL$ analogue of the 2-sections. 
Therefore, by defining their curvature, one defines the curvature of the $PL$ 
(polyhedral) 2-sections. Haantjes curvature (and, by passing to the dual complex 
whenever possible, Menger curvature as well) allows for what represents, we believe, 
a definition of sectional curvature having clear geometric content. This route not 
only allows an integrated and simple approach to the geometrization of networks, 
it is also, in a sense, the most natural and direct approach to this problem, 
given that $PL$ manifolds represent the geometric class closest to graphs or 
networks, on which proper and consistent differential geometry has been developed. 
Moreover, given the recent growing interest of the network community in the study 
of multiplex- and hypernetworks (as well as other higher dimensional generalizations 
of networks), viewing networks as $PL$ (and more general polyhedral manifolds) 
allows one not only to ration by analogy, but in fact provides us with arguably 
the best modeling tool. 

We would also like to add that given the lack of smoothness of the graph 
structure, hence of the full force of the differential geometric apparatus, 
each discretization of curvature captures, in the network context, only an 
essential aspect of the classical notion, and as such, it can satisfy only a 
restricted set of the full gamut of properties of the classical notion. Therefore 
it is not possible to find ``the best'' discrete curvature, but rather the best 
one suited for a certain task, in specific type of networks. This is exemplified 
in the sequel by the comparison between different notions of discrete Ricci 
curvature; see also the discussion on the comparative advantages of 
Ollivier-Ricci and Forman-Ricci curvature in \cite{Samal2018}. It is, therefore, 
desirable to devise different notions of discrete curvature and explore their 
specific advantages for various task. The present paper, represents, therefore, 
a step in this quest. 

The remainder of the paper is structured as follows. In Section 2, we introduce 
the Menger curvature. In Section 3, we introduce Haantjes curvature, its 
generalizations and extensions, as well as its use in the introduction of local 
Gauss-Bonnet theorem for networks. In Section 4, we provide a brief overview of 
the embedding properties of the metric curvatures presented here. In Section 5, 
we present empirical results from analysis of the Menger curvature and the 
Haantjes curvature in various model and real-world networks. In Section 6, we 
conclude with a summary and future outlook. In appendix \ref{sec:theory}, we 
present the mutual relationship between the two metric curvatures, the Menger 
curvature and the Haantjes curvature, and also their connections to the classical 
notion of curvature for curves. Some of the results reported in this manuscript 
were recently presented in a conference proceeding \cite{Saucan2019b}.


\begin{figure*}
\label{fig:spherical}
\centering
\includegraphics[width=5cm]{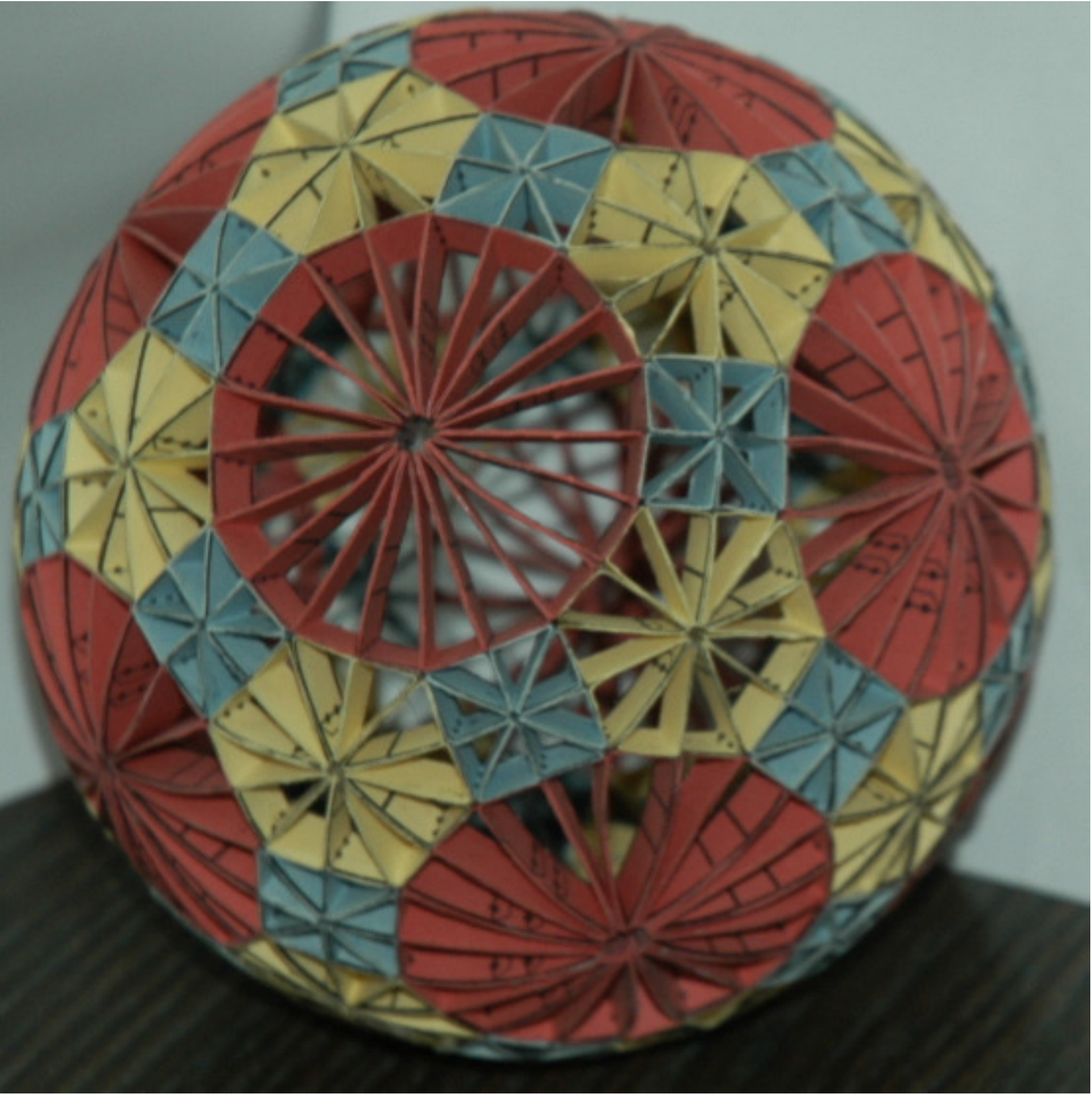}
\caption{The spherical network representing the 1-skeleton of the canonical 
triangulation with fundamental triangles of the spherical counterpart of the 
\textit{truncated icosidodecahedron} (see also Figure 3).}
\end{figure*}

\begin{figure*}
\label{fig:directedtriangle}
\centering
\includegraphics[width=10cm]{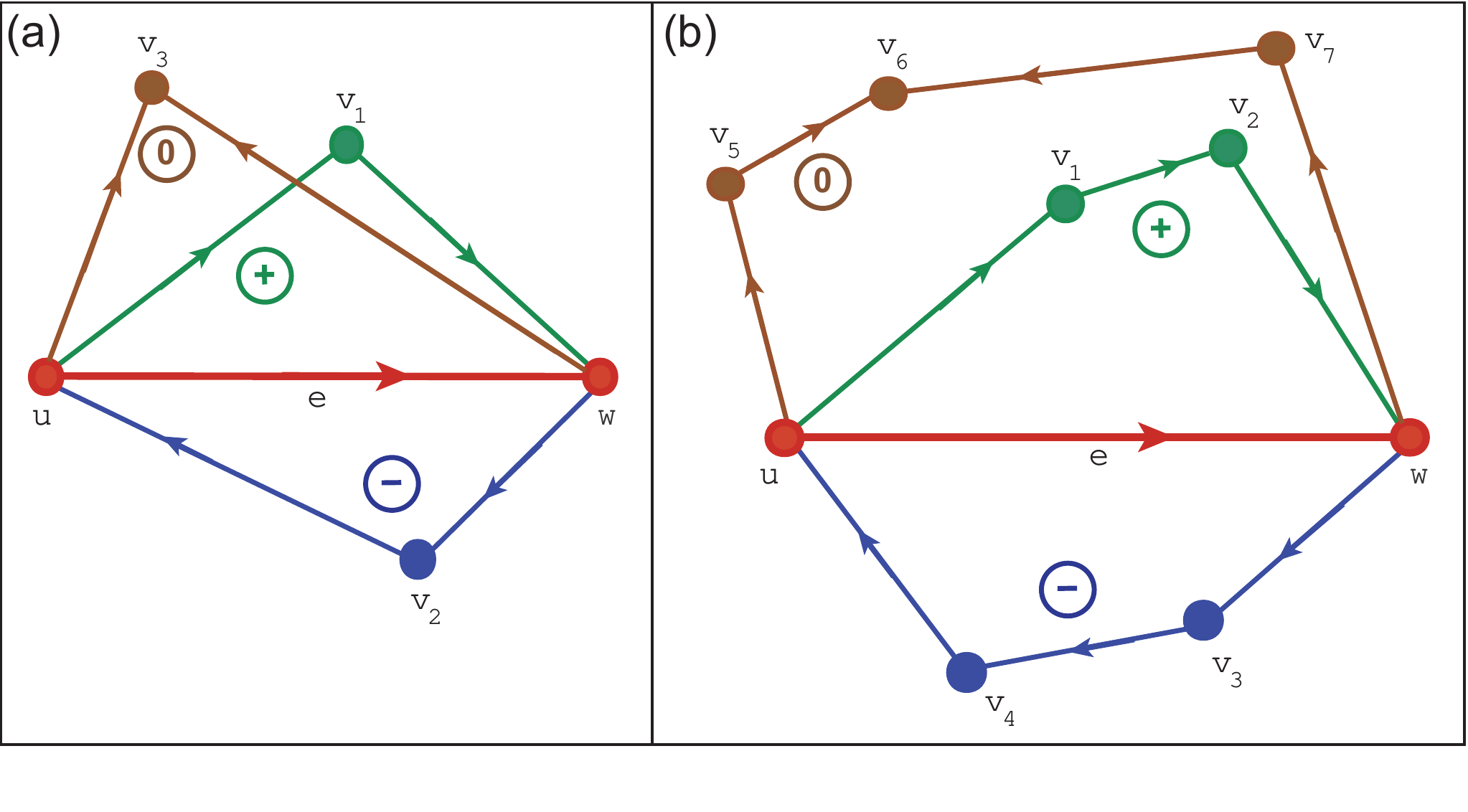}
\caption{Sign convention for {\bf (a)} directed triangles and {\bf (b)} directed 
polygons. This specific choice of the direction of triangles 
(while not the only possible one) is by now commonly accepted by many in the 
complex networks community (see e.g. \cite{Saucan2019b}), and it is motivated 
largely by the work of Alon and colleagues on network motifs \cite{Milo2002}. 
More precisely, positively oriented triangles correspond to feed forward loops, 
negatively oriented ones to feed backward loops, while zero curvature (no 
contribution) is attached to triangles that belong to neither of these types.
}
\end{figure*}

\section{Menger Curvature} 
\label{section:Menger}

The simplest, most elementary manner, of introducing curvature in metric spaces 
is due to Menger \cite{Menger1930}. Here, one simply defines the curvature 
$K(T)$ of a triangle $T$, i.e., a metric triple of points with sides of lengths 
$a, b, c$ as $\frac{1}{R(T)}$, where $R(T)$ is the radius of the circle 
circumscribed to the triangle. An elementary computation yields
\begin{equation} 
\label{eq:KM-Euc}
\kappa_{M,E}(T) = \frac{1}{R(T)} = \frac{abc}{4\sqrt{p(p-a)(p-b)(p-c)}}\,,
\end{equation}
where $p = \frac{(a+b+c)}{2}$ denotes the half-perimeter of $T$.

\vspace{0.25cm}
\noindent \textit{Remark:} The notion of curvature introduced above might appear 
somewhat counterintuitive, given that it prescribes to a triangle the curvature 
of its circumscribed circle. Thus, triangles of various shapes may have the same 
Menger curvature. Note that what appears as a more intuition of curvature, 
namely the so called {\it aspect ratio} is introduced in Section 3.3. While this 
notion is widely employed as such in Graphics (and, under other names, in Geometry 
and Geometric Analysis), it does not, however, represent a proper curvature notion 
(see Section 3.3 below).
\vspace{0.25cm}

However, there is conceptual problem with the above definition which utilizes 
the geometry of the Euclidean plane. In the general setting of networks, it is 
not natural to assume a Euclidean background. This is analogous to the geometry 
of surfaces, where the metric need not be Euclidean, but could be spherical, 
hyperbolic, or of varying Gauss curvature. For example, embedding networks in 
hyperbolic plane and space is becoming quite common \cite{Bianconi2017,Krioukov2010,
Zeng2010}.

Of course, one may formulate a spherical or hyperbolic analogue of Eq. 
\ref{eq:KM-Euc}, see e.g. \cite{Janson2015}. 
The spherical version is
\begin{equation} 
\label{eq:KM-Sph}
\kappa_{M,S}(T) = \frac{1}{\tan{R(T)}} = 
\frac{\sqrt{\sin{p}\sin{(p-a)}\sin{(p-b)}\sin{(p-c)}}}{2\sin{\frac{a}{2}}
\sin{\frac{b}{2}}\sin{\frac{c}{2}}}\,,
\end{equation}
whereas the hyperbolic version is
\begin{equation} 
\label{eq:KM-Hyp}
\kappa_{M,H}(T) = \frac{1}{\tanh{R(T)}} =
\frac{\sqrt{\sinh{p}\sinh{(p-a)}\sinh{(p-b)}\sinh{(p-c)}}}{2\sinh{\frac{a}{2}}
\sinh{\frac{b}{2}}\sinh{\frac{c}{2}}}\,.
\end{equation}

Note that, in the setting of networks, the constant factors ``4'' and ``2'',
respectively, appearing in the denominators of the above equations are less 
relevant, and they can be discarded in this context.

\vspace{0.25cm}
\noindent \textit{Remark:} The factors $(p-a), (p-b), (p-c)$ in the definition 
above are, in fact, the so called {\it Gromov products} of the vertices of $T$, 
which is employed in the definition of \textit{Gromov $\delta$-hyperbolicity} 
\cite{Gromov2007}. Thus Menger curvature (and a related, simpler metric 
invariant that we introduce in Section 3.3) may be viewed as a local version 
of the Gromov hyperbolicity. Given the fact that Gromov hyperbolicity is a 
global property, hence difficult to compute, it is good to be able to substitute, 
whenever possible, a  local related notion of curvature. Moreover, while Gromov 
hyperbolicity represents a coarse, global notion of sectional curvature, the 
Menger curvature which we propose herein is a local, notion of sectional curvature, 
thus fitting better the classical Riemannian paradigm of the development of 
curvature, first locally, then globally, that we adopted herein. This is especially 
meaningful if ``local-to-global'' properties, which play an important role in 
metric geometry \cite{Burago2001,Perelman1991,Plaut2002}, can be proven for 
Haantjes curvature.
\vspace{0.25cm}

\vspace{0.5cm}
\noindent \textbf{Examples of computing Menger curvature with spherical 
background geometry}
\begin{enumerate}
\item Let $\bigtriangleup ABC$ be a spherical triangle with angles $\alpha = 
\beta = \gamma = \frac{\pi}{2}$ on a sphere of radius 1. Then the length of 
the sides of the triangle, $a = b = c = \frac{\pi}{2}$, and the formula above 
renders $\kappa_{M,S}(\bigtriangleup ABC) = \frac{\sqrt{2}}{2}$.
\item The 1-skeleton of the triangulation of the sphere consisting of the 
fundamental triangles of a spherical counterpart of an Archimedean polyhedron 
(see Figure 1) represents a spherical network. The side of the spherical faces 
is taken to be 2 units. The fundamental triangle of an $n$-gonal face has 
angles $(\frac{\pi}{2},\frac{\pi}{3},\frac{\pi}{n})$, $n = 4,6,10$. Hence, 
the remaining sides of the triangles, thence their Menger curvatures, can be 
computed using the classical formulas of spherical trigonometry, see e.g. 
\cite{Janson2015}.
\end{enumerate}

\vspace{0.5cm}
\noindent \textbf{Examples of computing Menger curvature with hyperbolic 
background geometry.}
\begin{enumerate}
\item The hyperbolic plane admits a tessellation with triangles of angles 
$(\frac{\pi}{2},\frac{\pi}{3},\frac{\pi}{7})$. The dual graph represents also a 
tessellation with regular (hyperbolic) squares, hexagons and 14-gones, symbolically 
denoted as $(4,6,14)$. It can be viewed as the \textit{Cayley graph} of the group 
of the symmetries of the tessellation above, as generated by the reflections in 
the sides of the triangles, see e.g. \cite{Epstein1992}. Assuming that the edge is 
opposite the angle of measure $\frac{\pi}{7}$, one can compute the remaining sides 
using the classical formulas of hyperbolic trigonometry, see e.g. \cite{Janson2015}.
\item It is possible to tessellate the hyperbolic space with regular dodecahedra 
having all the faces with angles equal to $\frac{\pi}{2}$, see e.g. 
\cite{Thurston1997}. Then the fundamental triangle of the faces of these dodecahedra 
has angles $(\frac{\pi}{5},\frac{\pi}{4},\frac{\pi}{2})$. Normalizing the sides of 
the dodecahedra, such that half of its length be equal to 1, the \textit{hyperbolic 
law of sinuses}, see e.g. \cite{Janson2015}, gives $\sinh{a} = \frac{2\sqrt{2}}{2}
\cdot\frac{\sinh{1}}{\sinh{\pi/5}}$ for the side opposite the angle of measure 
$\frac{\pi}{4}$, and $\sinh{a} = \frac{\sinh{1}}{\sinh{\pi/5}}$ for the side opposite 
the right angle. From these the $\kappa_{M,H}(T)$ of the triangle can be readily 
computed.
\end{enumerate}

Hyperbolic geometry is considered better suited to represent the background 
network geometry as it captures the qualitative aspects of networks of exponential 
growth such as the World Wide Web, and thus, it is used as the setting for variety 
of purposes. However, spherical geometry is usually not considered as a model 
geometry for networks because that geometry has finite diameter, hence finite 
growth. However, spherical networks naturally arise in at least two instances. The 
first one is that of global communication, where the vertices represent relay 
stations, satellites, sensors or antennas that are distributed over the geo-sphere 
or over a thin spherical shell that can, and usually is, modeled as a sphere. The 
second one is that of brain networks, where the cortex neurons are envisioned, due 
to the spherical topology of the brain, as being distributed on a sphere or, in some 
cases, again on a very thin (only a few neurons deep) spherical shell, that can also 
be viewed as essentially spherical. One can also devise an analogous, although less 
explicit formula in spaces of variable curvature, but in network analysis, it is not 
clear where that background curvature should come from. After all, the purpose here 
is to define curvature, and not take it as given.

As defined, the Menger curvature is always positive. This may not be desirable, as 
in geometry, the distinction between positive and negative curvature is important. 
For directed networks, however, a sign $\varepsilon(T) \in \{-1,0,+1\}$ is naturally 
attached to a directed triangle $T$ (Figure 2), and the Menger curvature of the 
directed triangle is then defined, in a straightforward manner as
\begin{equation}
\kappa_{M,O}(T) = \varepsilon(T)\cdot\kappa_M(T)\,,
\end{equation}
where $\kappa_M$ could be the Euclidean, the spherical or the hyperbolic version, 
accordingly to the given setting. Note that, since in simplicial complexes, triangles 
adjacent to an edge represent discrete ($PL$) analogues of 2-dimensional sections, 
the Menger curvature of each such triangle $T$ can be naturally interpreted as the 
sectional curvature of $T$.

We can then define the Menger-Ricci curvature of an edge by averaging as in 
differential or piecewise linear geometry as
\begin{equation}
\kappa_{M,O}(e) = {\rm Ric}_{M,O}(e) = \sum_{T_e \sim e}\kappa_{M,O}(T_e)\,,
\end{equation}
where $T_e \sim e$ denote the triangles adjacent to the edge $e$, and the 
Menger-scalar curvature of a vertex is given by
\begin{equation} 
\label{eq:Menger-Ricci}
\kappa_{M,O}(v) = {\rm scal}_{M,O}(v) = 
\sum_{e_k \sim v} {\rm Ric}_{M,O}(e_k) = \sum_{T \sim v}{\kappa_{M,O}(T)}\,,
\end{equation}
where $e_k \sim v$ and $T \sim v$ stand for all the edges $e_k$ adjacent to the 
vertex $v$ and all the triangles $T$ having $v$ as a vertex, respectively.
Of course, for undirected networks, the sign of $\varepsilon(T)$ 
is assumed to be always equal to 1.

\vspace{0.25cm}
\noindent \textit{Remark:} ${\rm Ric}_M(e)$ captures, in keeping with the 
intuition behind $\kappa_M(T)$, the geodesic dispersion rate aspect of Ricci 
curvature. See \cite{Samal2018} for a succinct overview of the different aspects 
of Ricci curvature. 
\vspace{0.25cm}

As indicated above, we see two drawbacks for Menger curvature as a tool in network 
analysis. It depends on a background geometry model, and it naturally applies only 
to triangles, but not to more general 2-cells. Therefore, we turn, in the next section, 
to Haantjes  curvature, which we find more flexible in applications. (See, however, 
the remark below.) 

\vspace{0.25cm}
\noindent \textit{Remark:} It is possible to prescribe a Menger curvature for 
general paths as well, in the following manner. Let $\pi = v_0,v_1,\ldots,v_{n-1},
v_n$ be a path, subtended by the chord $(v_0,v_n)$, and let $v_k, 1 \leq k \leq n-1$ 
be any intermediary vertex. Then $v_k$ divides the path $\pi$ into two paths 
$\pi_1 = v_1,\ldots,v_k$ and $\pi_2 = v_k,\ldots,v_n$. Let us denote $a = l(\pi_1), 
b = l(\pi_2), c = l(v_0,v_n)$. Then one can consider the Menger curvature of the 
metric triangle $\bigtriangleup v_0v_kv_n$. Note that it does not depend on the 
choice of the vertex $v_k$. Again, the Menger curvature can be computed with the 
Euclidean, spherical or hyperbolic flavor, according to the preferred model geometry 
for the given network. As already noted above, this approach has the limitation of 
prescribing a predefined curvature for the network. On the other hand, this approach 
allows for the passing to an \textit{Alexandrov comparison} type of \textit{coarse 
geometry}, see e.g. \cite{Burago2001}.
\vspace{0.25cm}


\begin{figure*} 
\label{fig:archimedean}
\centering
\includegraphics[width=5cm]{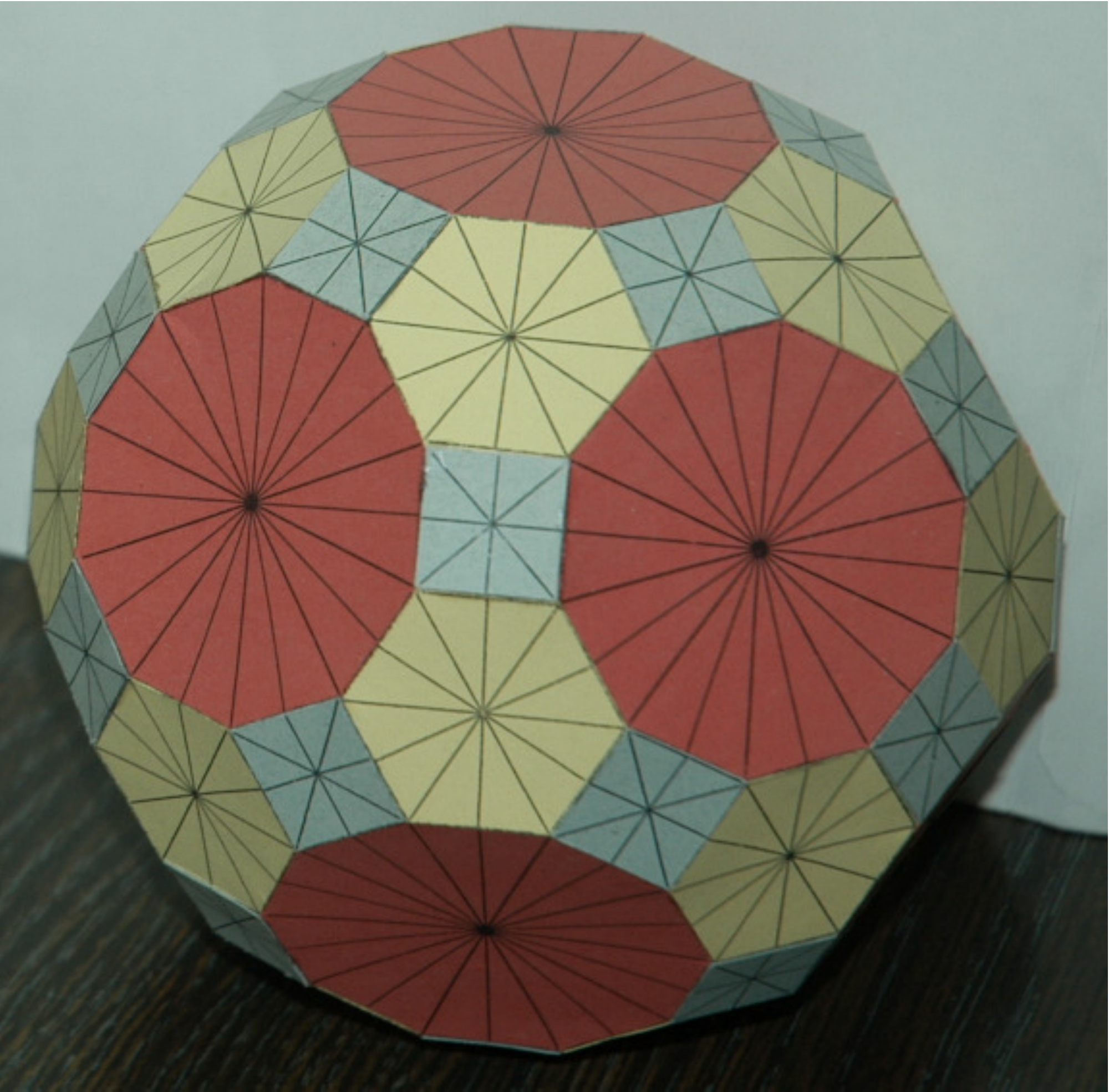}
\caption{The truncated icosidodecahedron, another Archimedean polyhedron. 
This is an Euclidean counterpart of the spherical polyhedron in Figure 1. 
There are three types of edges in this object, and the Forman-Ricci 
curvature for the square-hexagon, square-decagon and hexagon-decagon edges 
are equal to $\sqrt{2}+\sqrt{5}$, $\sqrt{2}+3$ and $\sqrt{5}+3$, 
respectively. Note that the computation of the Ollivier-Ricci 
curvature for this Archimedean polyhedron is not possible, due to the 
presence of faces with 6 or more edges, since the combinatorial 
Ollivier-Ricci curvature is not defined for such cycles. This fact 
underlines one of the important relative advantages of Haantjes-Ricci 
curvature, namely that it is applicable to networks containing cycles of 
any length.}
\end{figure*}
\begin{figure*}
\centering
\includegraphics[width=12cm]{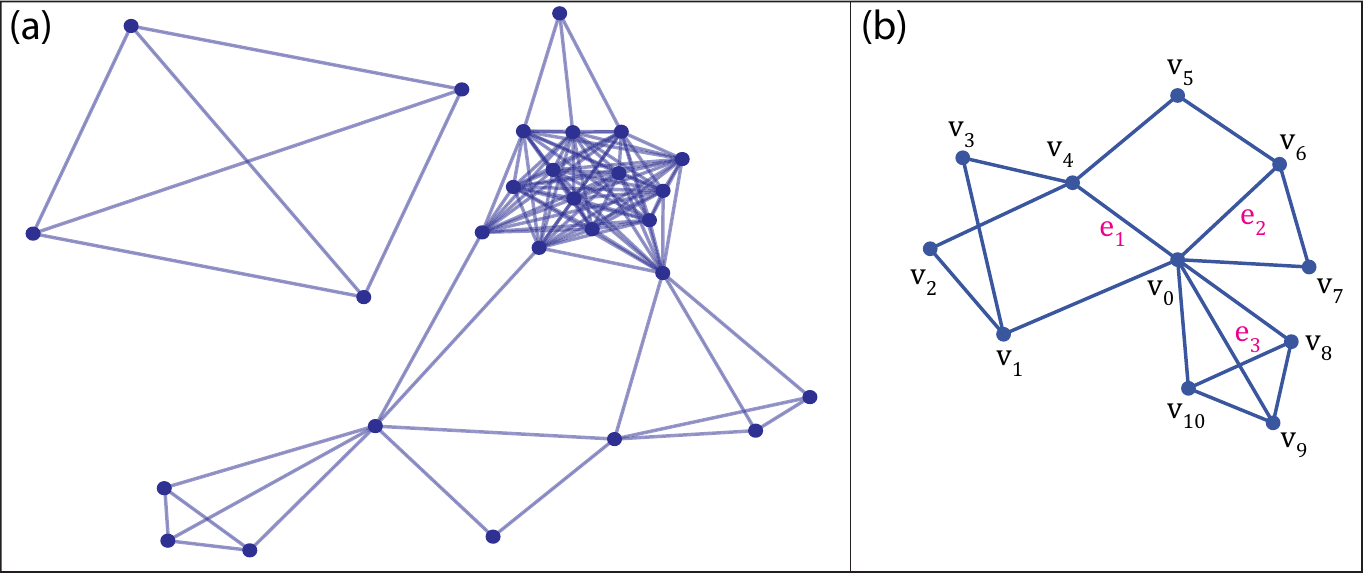}
\caption{An example of a classic social network Zebra 
\cite{Sundaresan2007}. {\bf (a)} Complete network. {\bf (b)} A 
portion of the Zebra network. Note that while for the edge $e_3$ 
both Menger- and Haantjes-Ricci curvature are computable, since the 
2-cycles adjacent to it are the triangles $\bigtriangleup v_0v_{10}v_9$ 
and $\bigtriangleup v_0v_{8}v_9$; for $e_2$ Menger-Ricci curvature can 
take into account only the triangle $\bigtriangleup v_0v_6v_7$, 
but not the quadrangle $\square v_0v_4v_5v_6$; while in the case of 
the edge $e_1$, Menger curvature would give a (false) 0 curvature, 
since there are no triangles, while Haantjes curvature ``sees'' the 
3 adjacent quadruples, namely $\square v_0v_1v_2v_4$, $\square v_0v_1v_3v_6$, 
$\square v_0v_4v_5v_6$. This simple case exemplifies, therefore, the greater 
flexibility of Haantjes curvature, that is not restricted, as Menger's, to 
simplicial complexes.}	
\end{figure*}

\begin{figure*}
\centering
\includegraphics[width=12cm]{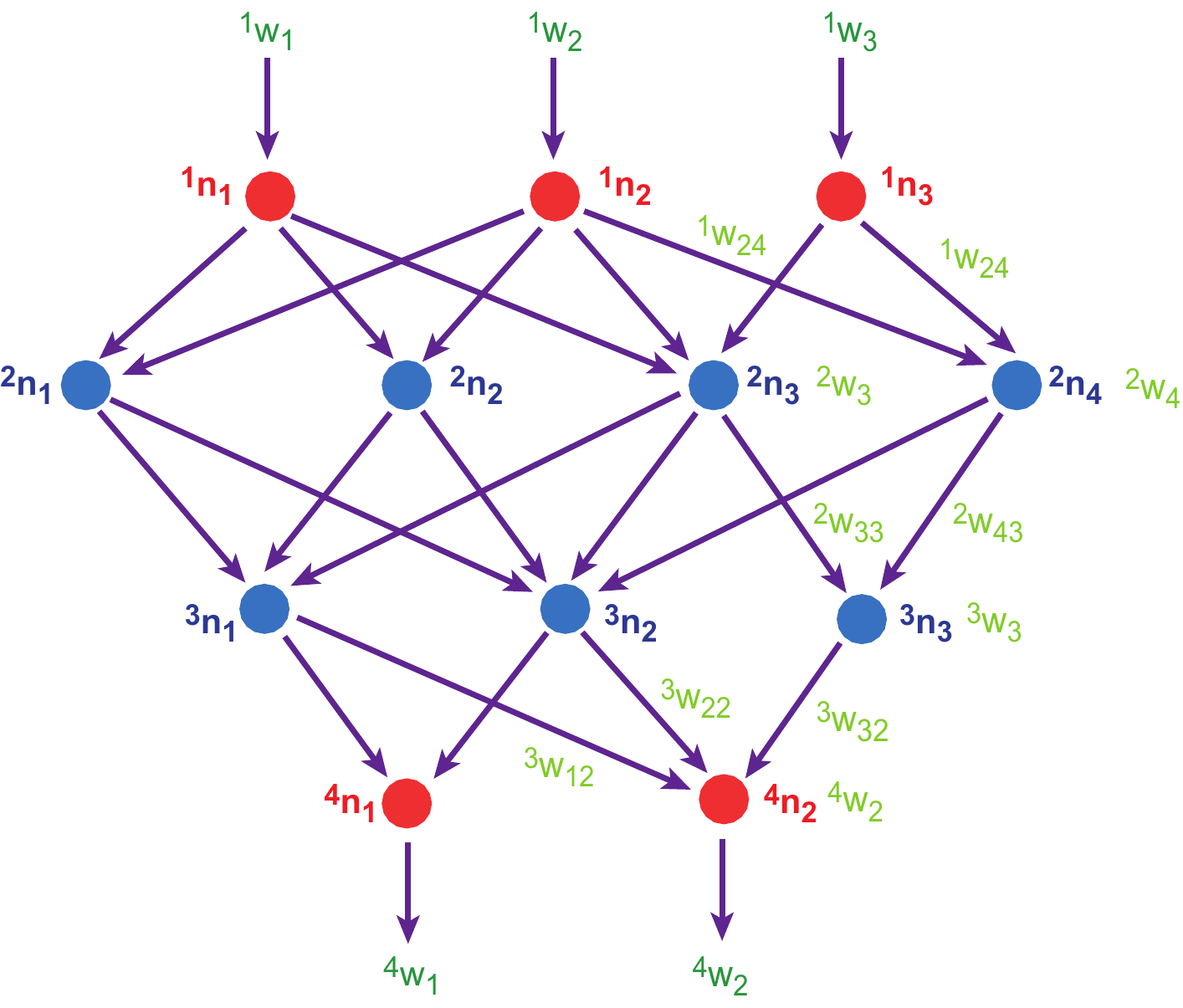}
\caption{Schematic illustration of an artificial neural network 
(ANN). The nodes of the input layer, namely ${^1}n_1$, ${^1}n_2$, ${^1}n_3$  
have attached weights, ${^1}w_1$, ${^1}w_2$, ${^1}w_3$ $\geq 0$, representing 
the respective inputs at these nodes. Each of the synapses ${^i}s_{jk}$, i.e. 
directed edges, has a weight ${^i}w_{jk} \geq 0$, and these weights, together 
with the transfer functions ${^i}f_k$, determine the weights ${^i}w_{jk} \geq 
0$, of the neurons ${^i}n_k$, up to and including the weights ${^4}w_1, {^4}w_2$ 
of the output layer ${^4}n_1, {^4}n_2$. Observe that there are no proper 
triangles, thus Forman-Ricci curvature can not be computed for beyond separate 
edges and, given the essentially tree-like structure of the network, it is 
not very expressive; nor is the Menger-Ricci curvature definable. However, 
the Haantjes-Ricci curvature in any direction (of length 2) is definable. 
Just to exemplify, if we consider the synapses ${^2}s_{11}$ and ${^2}s_{12}$ 
to be endowed with combinatorial weight 1, and the weights ${^3}w_{12}, {^3}w_{22}$ 
to be equal to 2 and 3, respectively, the Haantjes-Ricci curvature in the 
direction ${^2}n_1 {^3}n_1 {^4}n_2$ will be $13\sqrt{3}$. Needles to say, 
in order to truly understand the significance of this type of curvature, 
one need to systematically explore the various types of ANN architectures 
as well as the diverse transfer functions specific to each learning task, 
a goal that clearly goes beyond the scope of the present article.}	
\end{figure*}
\begin{figure*}
\label{fig:haantjes-path}
\centering
\includegraphics[width=5.5cm]{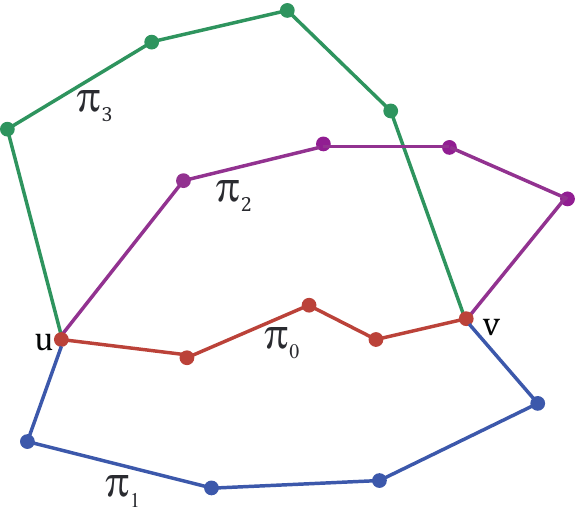}
\caption{Haantjes-Ricci curvature in the direction $\overline{uv}$, is defined as 
${\rm Ric}_H(\overline{uv}) = \sum_1^mK^i_{H,\overline{uv}} = 
\sum_1^m\kappa_{H,\overline{uv}}(\pi_i)$\, where $\kappa^{2}_{H,\overline{uv}}(\pi_i) 
= \frac{l(\pi_i) - l(\pi_{0})}{ l(\pi_{0})^3}\, i = 1,2,3$; $\pi_0$ being the 
shortest path connecting the vertices $u$ and $v$.}
\end{figure*}

\begin{figure*}
\label{fig:example-haantjes}
\centering
\includegraphics[width=11cm]{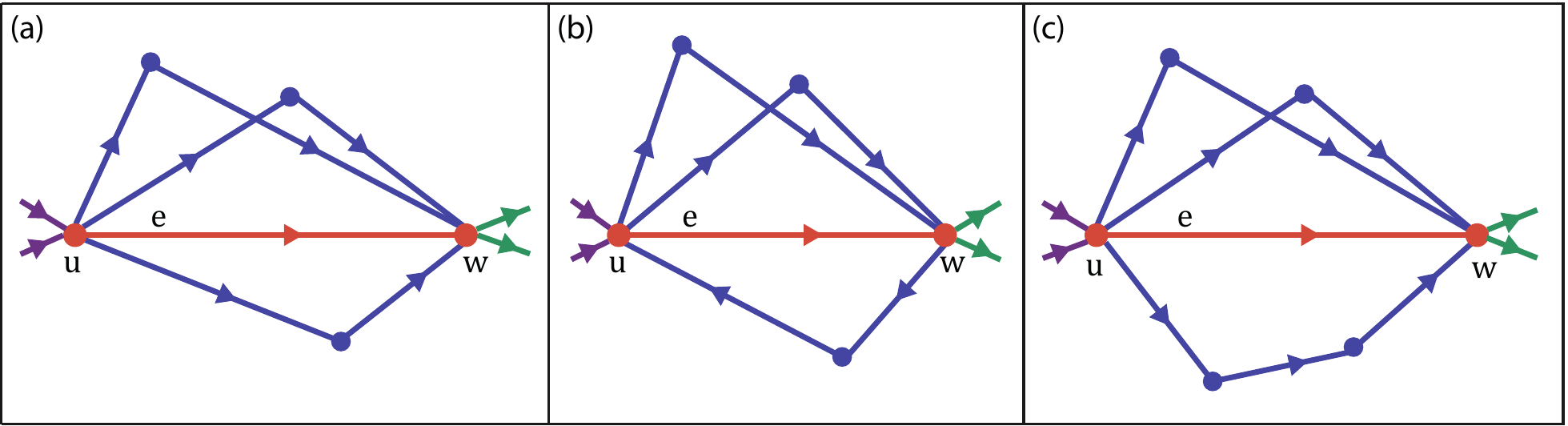}
\caption{Haantjes-Ricci curvature  ${\rm Ric}_{H,O}$ for directed networks, endowed 
with combinatorial weights: {\bf (a)} ${\rm Ric}_{H,O} = 6\pi - 3$, {\bf (b)} 
${\rm Ric}_{H,O} = 2\pi - 1$, and {\bf (c)} ${\rm Ric}_{H,O} = 6\pi - 2 - \sqrt{2}$. 
Note that while Forman-Ricci curvature is a counter of triangles in simplicial 
complexes, Haantjes-Ricci curvature represents a counter of all $n$-gones, 
since each $n$-gone contributes a $\sqrt{n}$ term.}
\end{figure*}

\section{Haantjes Curvature} 
\label{section:Haantjes}

Haantjes \cite{Haantjes1947} defined metric curvature by comparing the ratio 
between the length of an arc of curve and that of the chord it subtends. More 
precisely, if $c$ is a curve in a metric space $(X,d)$, and $p, q, r$ 
are points on $c$, $p$ between $q$ and $r$, the Haantjes curvature is defined as
\begin{equation}                         
\label{eq:haantjes-1}
\kappa_{H}^2(p) = 24\lim_{q,r \rightarrow p}\frac{l(\widehat{qr})-d(q,r)}
{\big(d(q,r)\big)^3}\,,
\end{equation}
where $l(\widehat{qr})$ denotes the length, in the intrinsic metric induced by 
$d$, of the arc $\widehat{qr}$ (see also Appendix below). In the network case, 
$\widehat{qr}$ is replaced by a path $\pi = v_0,v_1,\ldots,v_n$, and the 
subtending chord by edge $e = (v_0,v_n)$. Clearly, the limiting process has no 
meaning in this discrete case. Furthermore, the normalizing constant ``24'' which 
ensures that, the limit in the case of smooth planar curves will coincide with 
the classical notion, is superfluous in this setting. This leads to the following 
definition of the \textit{Haantjes curvature of a simple path} $\pi$:
\begin{equation}                         
\label{eq:haantjes-path}
\kappa_{H}^2(\pi) = \frac{l(\pi)-l(v_0,v_n)}{l(v_0,v_n)^3}\,,
\end{equation}
where, if the graph is a metric graph, $l(v_0,v_n) = d(v_0,v_n)$. In particular, 
in the case of the combinatorial metric, we obtain that, for path $\pi = v_0,v_1,
\ldots,v_n$ as above, $\kappa_H(\pi) = \sqrt{n-1}$. Note that considering simple 
paths is not a restriction, given that a metric arc is, by definition, a simple 
curve. However, to capture in the discrete context the local nature of the Ricci 
(and scalar) curvature, we shall restrict to paths $\pi$ such that $\pi^{*} = v_0,
v_1,\ldots,v_n,v_0$ is an {\it elementary cycle}.

Clearly, one can extend the above definition of Haantjes curvature to directed 
paths in the same manner as done for the Menger curvature in the preceding section, 
namely
\begin{equation}
\kappa_{H,O}(\pi) = \varepsilon(\pi)\cdot\kappa_H(\pi)\,,
\end{equation}
for every directed path $\pi$, where $\varepsilon \in \{-1,0,+1\}$ denotes the 
direction of path $\pi$ (see also Section 3.3 below). As in the 
case of Menger curvature, $\varepsilon(\pi) \equiv +1$, for undirected networks.

In a straightforward manner, we can define Haantjes-Ricci curvature and 
Haantjes-scalar curvature, similar to the Menger curvature, as
\begin{equation}
\label{eq:haantjes-compute}
\kappa_{H,O}(e) = {\rm Ric}_{H,O}(e) = \sum_{\pi \sim e}\kappa_{H,O}(\pi)\,,
\end{equation}
where $\pi \sim e$ denote the paths that connect the vertices anchoring the edge 
$e$, and
\begin{equation}
\kappa_{H,O}(v) = {\rm scal}_{H,O}(v) = \sum_{e_k \sim v} {\rm Ric}_{H,O}(e_k)\,,
\end{equation}
where $e_k \sim v$ stands for all the edges $e_k$ adjacent to the vertex $v$. In the 
following, we shall however significantly strengthen the somewhat simplistic definition 
of the Haantjes-Ricci curvature presented above.

\vspace{0.5cm}
\noindent \textbf{Examples of Archimedean (semi-regular) polyhedra.}
\begin{enumerate}
\item \textit{The truncated dodecahedron} is a typical Archimedean convex body, 
hence of positive combinatorial curvature concentrated at vertices. However, the 
two types of edges, the pentagon-hexagon ones and the hexagon-hexagon ones, in 
this object have Forman-Ricci curvature ${\rm Ric}_F(e)$ equal to -1 and -2, 
respectively. In contrast, the Haantjes-Ricci curvature is, by definition, always 
positive: ${\rm Ric}_H(e)=9$ for the pentagon-hexagon edges, and ${\rm Ric}_H(e)=10$, 
for the hexagon-hexagon edges.
\item \textit{The truncated octahedron} is another Archimedean polyhedron. The 
positivity condition for Forman-Ricci curvature for all edges does not hold here 
as well: ${\rm Ric}_F(e)=0$ for square-hexagon edges, and ${\rm Ric}_F(e)=-2<0$ 
for hexagon-hexagon edges. While the Haantjes-Ricci curvature is, by definition, 
always positive: ${\rm Ric}_H(e)=8$ for the square-hexagon edges, and ${\rm Ric}_H(e)=10$,
for the hexagon-hexagon edges. Let us also note that the edges of this object represent 
the \textit{Cayley graph} of the symmetric group $S(4)$ with respect to the set of 
generators $\{\tau_1,\tau_2,\tau_3\}$, where $\tau_1,\tau_2,\tau_3$ are the 
transpositions $\tau_1 = (12), \tau_2 = (23), \tau_3 = (34)$ (see \cite{Epstein1992}).
\end{enumerate}

It should be noted that the computation of the Ollivier-Ricci 
curvature for the Archimedean polyhedra above is not always possible, as they 
contain faces with 6 or more edges, and the combinatorial Ollivier-Ricci 
curvature is not defined for such cycles (Figure 3). (See also Table 1, where 
the Ollivier-Ricci curvature of the hexagonal planar grid is omitted for the 
very same reason.) This fact demonstrates one of the important relative advantages 
of Haantjes-Ricci curvature, namely that it is applicable to (networks 
containing) cycles of any length. This feature of Haantjes-Ricci curvature 
is of importance not only in a geometric context, such in the examples above, 
but also in such settings as that of social and biological networks, for 
instance, where the length of cycles has important modeling significance.

\vspace{0.5cm}
\noindent \textbf{Examples of non-convex uniform polyhedra.}
\begin{enumerate}
\item \textit{The tetrahemihexahedron} $(4.\frac{3}{2}. 4. 3)$ is a simple 
non-orientable polyhedron, representing a (minimalistic) model of the real 
projective plane. As such, it's Euler characteristic is equal to 1. The 
Forman-Ricci curvature of each edge here is equal to 1, and thus, strictly 
positive, whereas the Haantjes-Ricci curvature is equal to $1+\sqrt{2}$.
\item \textit{The octahemoctahedron} $(6.\frac{3}{2}.6.3)$ is the only 
toroidal uniform poyhedron with Euler characteristic 0. The Forman-Ricci 
curvature of each triangle-hexahedron edge here is equal to -1, whereas the 
Haantjes-Ricci curvature is equal to $3$. (Note that for this 
uniform polyhedron the Ollivier-Ricci curvature is, again, not defined, due 
to the existence faces with 6 edges.)
\item \textit{The nonconvex great rhombicuboctahedron (or 
quasi-rhombicuboctahedron)} $(4.\frac{3}{2}.4.4)$ is composed of 8 regular 
triangles and 18 squares. The triangles are \textit{retrograde}, and thus, 
their Haantjes-Ricci curvature should be taken with a \textit{negative} sign. 
Therefore, the Haantjes-Ricci curvature of an triangle-square edge here 
is $\sqrt{2}-1$, and not $\sqrt{2}+1$.
\end{enumerate}

\vspace{0.5cm}
\noindent \textbf{Examples of three-dimensional polyhedral complexes.}
\begin{enumerate}
\item \textit{The Seifert-Weber dodecahedral space} is a three manifold obtained 
by gluing, via $\frac{3}{10}$ of a clockwise full twist, the opposite faces of a 
regular dodecahedron. Each edge is incident to 5 dodecahedra/pentagonal faces, and 
hence, the Forman-Ricci curvature of such edges is $-7$, thus corresponding to the 
fact that it is possible to tile the Hyperbolic space with such dodecahedra, 
whereas the Haantjes-Ricci curvature is equal to ${\rm Ric}_H(e)=5\sqrt{3}$.
\item \textit{The Poincar\'{e} dodecahedral space} is a three manifold obtained 
by gluing, via $\frac{1}{10}$ of a clockwise full twist, the opposite faces of a 
regular dodecahedron. Each edge is incident to 3 dodecahedra/2-faces, and hence, 
the Forman-Ricci curvature of such edges is equal to $-1$, even though the 
Poincar\'{e} is the classical homology 3-sphere, and thus, expected to have a 
positive curvature. On the other hand, the Haantjes-Ricci curvature of such edges 
is positive, more precisely, ${\rm Ric}_H(e)=3\sqrt{3}$.
\end{enumerate}
\vspace{0.5cm}
\noindent \textbf{Example of a classic social network.} A small, classic 
example of a social network, Zebras \cite{Sundaresan2007}, exemplifies 
the fact that, given that it is defined only on triangles, Menger-Ricci
curvature is not always computable, or it gives only very partial, hence 
distorted results (Figure 4). However, Haantjes-Ricci curvature can be 
easily computed for this network.

\vspace{0.5cm}
\noindent \textbf{Example of a Artificial Neural Network (ANN)} Figure 5 shows 
a small example of a weighted ANN which demonstrates the advantages of 
Haantjes-Ricci curvature, over Menger-Ricci, for the intelligence of networks 
arising in the deep learning context.

\vspace{0.25cm}
\noindent \textit{Remark:} Among metric curvatures, the Haantjes curvature has 
several advantanges from the perspective of network science applications. Aside 
from the simplicity of computation in networks, Haantjes-Ricci curvature for an 
edge in undirected and unweighted (combinatorial) simplicial complexes can be 
obtained by simply counting the triangles $t$ containing an edge $e$, that is   
\begin{equation}
{\rm Ric}_H(e) = \sharp\{t\,|\, t > e\}\,.
\end{equation}
Moreover, in $k$-regular undirected and unweighted simplicial complexes where 
each vertex is incident to precisely $k$ edges, the (augmented) Forman-Ricci 
curvature of an edge $e$ \cite{Samal2018} is given by	
\begin{equation}
{\rm Ric}_F(e) = 4 - 3\sharp\{t\,|\, t > e\} - 2k\,,
\end{equation}
and thus, the formula also reduces to counting the triangles $t$ containing an 
edge $e$. While the later formula above is also dependent on vertex degree $k$, 
for any given  $k$ the two types of Ricci curvature, Haantjes and Forman, for 
edges in undirected combinatorial simplicial complexes both reduce to the 
counting of triangles adjacent to a given edge. The computation of Haantjes-Ricci 
curvature above has an additional advantage of not being dependent on $k$. On the 
other hand, note that Haantjes-Ricci curvature is always positive while 
Forman-Ricci curvature is mostly negative. Still, the distributions of the two 
curvatures, over a network, are likely to strongly correlated. We should, however, 
note here that ${\rm Ric}_F$ covers further aspects of network geometry that 
${\rm Ric}_H$ does not. (However, see also the discussion in the concluding 
section.)
\vspace{0.25cm}

\begin{table*}
\label{table:comparisons}
\caption{Comparison of undirected curvatures for a number of standard grids 
(tessellations) of the Euclidean plane and space.}
\begin{tabular}{l@{\hskip 0.5in}c@{\hskip 0.2in}c@{\hskip 0.2in}c@{\hskip 0.2in}c}
\hline \hline
Curvature	& Triangular	& Square  		& Hexagonal 	& Euclidean \\
Type 		& Tessellation	& Tessellation 	& Tessellation	& Cubulation \\ 
\hline
${\rm Ric}_H(e)$ & {4$\pi$ - 2} & {4$\pi$ - 4} & {4$\pi - 2\sqrt{2}$} & {8$\pi - 4\sqrt{2}$} \\ 
${\rm Ric}_{F,r}(e)$ & {-8} & {-2} & {-2} & {-4} \\ 
${\rm Ric}_{F}(e)$ & {-2} & {0} & {4} & {4} \\ 
${\rm Ric}_{O}(e)$ & {1} & {-1} & {--} & {-$\frac{4}{3}$}\\ 
\hline \hline
\end{tabular}	
\end{table*}	

\subsection{Local Gauss-Bonnet Theorem based sectional and  Ricci curvatures}

Due to its advantages over Menger curvature, we shall now use Haantjes 
curvature to provide stronger definitions of scalar curvature and Ricci 
curvature of networks. Here, the basic idea is to adapt the \textit{local 
Gauss-Bonnet theorem} to this discrete setting. This approach 
allows us to further advance our program of providing a $PL$ type geometry 
for networks and their higher dimensional counterparts.
Recall that, in the classical context of smooth surfaces, the theorem 
states that
\begin{equation}  
\label{eq:SmoothGB}
\int_DKdA + \sum_0^p\int_{v_i}^{v_{i+1}}k_gdl + \sum_{0}^p\varphi_i 
= 2\pi\chi(D)\,,
\end{equation}
where $D \simeq \mathbb{B}^2$ is a (simple) region in the surface, having as 
boundary $\partial D$ a piecewise-smooth curve $\pi$, of vertices (i.e., points 
where $\partial D$ is not smooth) $v_i, i=1,\dots,n$, ($v_n = v_0$), 
$\varphi_i$ denotes the external angles of $\partial D$ at the vertex $v_i$, 
and, $K$ and $k_g$ denote the Gaussian and geodesic curvatures, respectively.

Let us first note that, in the absence of a background curvature, the very notion 
of angle is undefinable. Thus, for abstract (non-embedded) cells, there exists no 
\textit{honest} notion of angle. Hence, the last term on the left side of above 
Eq. \ref{eq:SmoothGB} has no proper meaning, and thus, can be discarded. Moreover, 
the distances between non-adjacent vertices on the same cycle (apart from the path 
metric) are not defined, and thus, the second term on the left side of above Eq. 
\ref{eq:SmoothGB} also vanishes.

We next concentrate on the case of combinatorial (unweighted) networks. For such  
networks endowed with the combinatorial metric, the area of each cell is usually 
taken to be equal to 1. Moreover, one can naturally assume that the curvature is 
constant on each cell, and thus, the first term on the left side of Eq. 
\ref{eq:SmoothGB} reduces simply to $K$. In addition, given that $D$ is a 2-cell, 
we have $\chi(D) = 1$. Therefore, in the absence of the definition of an angle, 
it is naturally to define 
\begin{equation}
K = 2\pi - \int_{\partial D}k_gdl\,.
\end{equation}
It is tempting to next consider $\partial D$ as being composed of segments (on 
which $k_g$ vanishes), except at the vertices, thus rendering the above expression 
as
\begin{equation}
\label{eq:haantjes-5}
K = 2\pi - \sum_{1}^n\kappa_H(v_i)\,.
\end{equation}

\vspace*{0.25cm}
\noindent {\it Remark}
An alternative approach to defining the curvature of cell would be the following: 
Since in a Euclidean polygon, the sum of the angles equals $\pi (n-2)$, where $n$ 
represents the number of vertices of the polygon, one could replace the angle sum 
term in Eq. \ref{eq:SmoothGB} simply by $\pi (n-2)$.
\vspace*{0.25cm}

Moving to the general case of weighted networks, one can not define a (non-trivial) 
Haantjes curvature for vertices, as already noted above, no proper distance between 
two non-adjacent vertices $v_{i-1}$ and $v_{i+1}$ on the same cycle can be 
implicitly assumed (apart from the one given by the path metric, which would 
produce trivial zero curvature at vertex $v_i$). In fact, in this general case, 
neither can the arc (path) $\pi = v_0,v_1,\ldots,v_n$ be truly viewed as smooth. 
Therefore, we have no choice but to replace the second term on the right side of 
above Eq. \ref{eq:haantjes-5} by $\kappa_H(\pi)$, where it should be remembered 
that $\pi$ represents the path $v_0,v_1,\ldots,v_n$ of chord $e = (v_0,v_n)$. 

We can now define the \textit{Haantjes-sectional curvature of a 2-cell} 
$\mathfrak{c}$. Given an edge $e = (u,v)$ and a 2-cell $\mathfrak{c}$, $\partial 
\mathfrak{c} = (u=v_0,v_1,\ldots,v_n=v)$ (\textit{relative to the edge $e \in \partial 
\mathfrak{c}$}), we have
\begin{equation} 
\label{eq:K=-kH}
K_{H,e}(\mathfrak{c}) = 2\pi -\kappa_{H,e}(\pi)\,,
\end{equation}
where $\pi$ denotes the path $v_0,v_1,\ldots,v_n$ subtended by the chord 
$e = (v_0,v_n)$ (see Figure 6), and $\kappa_{H,e}(\pi)$ denotes its 
respective Haantjes curvature. In the sequel we shall refer 
to this version as the {\it strong} Haantjes-Ricci curvature of the cell 
$\mathfrak{c}$ to distinguish it from the simpler version each time 
such a differentiation is required.

Note that the definition above is much more general than the one based on Menger 
curvature. Indeed, not only is it applicable to cells whose boundary has (combinatorial) 
length greater than three, i.e., not just to triangles, it also does not presume any 
convexity condition for the cells, even in the case when they are realized in some model 
space, e.g. in $\mathbb{R}^3$. However, for simplicial complexes endowed with the 
combinatorial metric, the two notions coincide up to a constant. More precisely, in this 
case, for any triangle $T$, $\kappa_M(T) \slash \kappa_H(T) = \sqrt{3}/3$. In fact, for
the case of smooth, planar curves, Menger and (unnormalized) Haantjes curvature coincide 
in the limit, and furthermore, they agree with the classical concept. However, for 
networks there is no proper notion of convergence, a fact which allowed us to discard 
the factor 24 in the original definition (Eq. \ref{eq:haantjes-1}) of Haantjes curvature.

We can now define, analogous to Eq. \ref{eq:Menger-Ricci}, the \textit{Haantjes-Ricci 
curvature of an edge} $e$ as
\begin{equation} 
\label{eq:Haantjes-Ricci}
{\rm Ric}_H(e) = \sum_{\mathfrak{c} \sim e}K_{H,e}(\mathfrak{c}) =
\sum_{\mathfrak{c} \sim e}(2\pi -\kappa_{H,e}(\pi))\,,
\end{equation}
where the sum is taken over all the 2-cells $\mathfrak{c}$ adjacent to $e$.
Again, we shall call this version  the {\it strong} Haantjes-Ricci 
curvature of the edge $e$ to differentiate it from the first, simpler version, 
in each instance where the confusion is possible.
See Figure 7 for examples of computation of Haantjes-Ricci curvature in networks, 
in the directed case. See Table 1 
\footnote{Note that the field assigned to the hexagonal tilling for Ollivier-Ricci 
curvature in Table 1 is marked as ``--'', since in this case the Ollivier-Ricci 
curvature is not applicable, see \cite{Cushing2019}. Also, in 
Table 1, ${\rm Ric}_H(e)$ denotes the strong Haantjes-Ricci curvature given by 
Eq. \ref{eq:Haantjes-Ricci}.}
for the comparison on a number of (undirected) standard planar and spatial grids 
of the various types of Ricci curvature at our disposal.


\subsection{The case of general weights}

In this subsection, we return to the general case of weighted graphs. Firstly, 
note that it is not reasonable to attach area 1 to every 2-cell in such graphs. 
However, as discussed in \cite{Horak2013,Saucan2018}, it is possible to endow 
cells in an abstract weighted graph with weights that are both derived from the 
original ones and have a geometric content. For instance, in the case of 
unweighted social or biological networks, endowed with the combinatorial metric,
one can designate to each face, instead of the canonical combinatorial weight 
equal to 1, a weight that \textit{penalizes} the faces with more edges, and thus, 
reflecting the weaker mutual connections between the vertices of such a face.
Thus, it is possible to derive a proper local Gauss-Bonnet formula for such 
general networks, in a manner that still retains the given data, yet captures the 
geometric meaning of area, volume, etc. Thus, when considering any such geometric 
weight $w_g(\mathfrak{c}^2)$ of a 2-cell $\mathfrak{c}^2$, the appropriate form 
of the first term on the left side of Eq. \ref{eq:SmoothGB} becomes
\begin{equation}
\nonumber
Kw_g(\mathfrak{c}^2)\,,
\end{equation}
and the fitting form of Eq. \ref{eq:K=-kH}
\begin{equation} 
\label{eq:K-fitting}
K_{H,e}(\mathfrak{c}) = -\frac{1}{w_g(\mathfrak{c}^2)}\left(2\pi - \kappa_{H,e}(\pi)\right)\,.
\end{equation}
Before passing to the problem of extending the above definition to the case of 
general weights, let us note that the observations above regarding Menger curvature 
for directed networks apply also to Haantjes curvature, after properly extending 
the notion to directed 1-cycles of any length and not just to directed triangles 
(see Figure 2). Again, as for the Menger curvature, considering directed networks 
actually simplifies the problem, in the sense that it allows for variable curvature 
(and not just one with constant sign). For general edge weights, we have the problem 
that the \textit{total weight} $w(\pi)$ of a path $\pi = v_0,v_1,\ldots,v_n$ is not 
necessarily smaller than the weight of its subtending chord $e=(v_0,v_n)$ 
\footnote{We suggest the name \textit{strong local metrics} for those sets of 
positive weights that satisfy the generalized triangle inequality 
$w(v_0,v_1,\ldots,v_n,v_{n+1}) < w(v_0,v_n)$, 
for any elementary 1-cycle $v_0,v_1,\ldots,v_n,v_0$\,.}, 
thus Haantjes' definition cannot be applied. However, we can turn this to our own 
advantage by reversing the roles of $w(\pi)$ and $w(v_0,v_n)$ in the definition of
the Haantjes curvature and assigning a minus sign to the curvature of cycles for 
which this occurs. Thus, this approach actually allows us to define a variable sign 
Haantjes curvature of cycles (hence, a Ricci curvature as well), even if the given 
network is not a naturally directed one.

Note that the case when $w(v_0,v_1,\ldots,v_n) = w(v_0,v_n)$, i.e., that of zero 
curvature of the 2-cell $\mathfrak{c}$ with $\partial \mathfrak{c} = v_0,v_1,\ldots,v_n,v_0$ 
straightforwardly corresponds to the splitting case for the path metric induced by 
the weights $w(v_i,v_{i+1})$.

Indeed, the method suggested above reduces to the use of the path metric, in most of 
the cases. One can always pass to the path metric and apply to it the Haantjes curvature. 
Beyond the complications that this might induce in certain cases, it is, in our view, 
less general, at least from a theoretical viewpoint, since it necessitates the passage 
to a metric. However, in the case of most general weights, i.e. both vertex and edge 
weights, one has to pass to a metric. We find the \textit{path degree metric} (see e.g. 
\cite{Keller2015}) especially alluring as it is both simple and has the capacity to 
capture, in the discrete context, essential geometric properties of Riemannian metrics. 
However, we also refer the reader to \cite{Saucan2005} for an ad hoc metric devised 
precisely for use on graphs in tandem with Haantjes curvature.


\subsection{A further generalization}

Note that Eqs. \ref{eq:haantjes-1} and \ref{eq:haantjes-path} are meaningful not only 
for a single edge, we can consider any two vertices $u,v$ that can be connected by a 
path $\pi$. Among the simple paths $\pi_1,\ldots,\pi_m$ connecting the vertices, the 
shortest one, i.e., the one for which $l(\pi_{i_0}) = \min\{l(\pi_1),\ldots,l(\pi_m)\}$ 
is attended represents the \textit{metric segment} of ends $u$ and $v$. Therefore, 
given any two such vertices, we can define the Haantjes-Ricci curvature in the direction 
$\overline{uv}$ to be
\begin{equation} 
\label{eq:haantjes-generalized}
{\rm Ric}_H(\overline{uv}) =  \sum_1^mK^i_{H,\overline{uv}} 
= \sum_1^m\kappa_{H,\overline{uv}} (\pi_i)\,
\end{equation}
where $K^i_{H,\overline{uv}}$ denotes the Haantjes-Ricci curvature of the cell 
$\mathfrak{c}_i$, where $\partial \mathfrak{c}_i = \pi_i\pi_0^{-1}$, relative to the 
direction $\overline{uv}$, and where
\begin{equation}
\kappa^2_H(\pi_i) = \frac{l(\pi_i) - l(\pi_{0})}{ l(\pi_{0})^3}\,,
\end{equation}
and where the paths $\pi_1,\ldots,\pi_m$ satisfy the condition that $\pi_i\pi_0^{-1}$ 
is an \textit{elementary cycle}. This represents a locality condition in the network 
setting.

\vspace{0.25cm}
\noindent \textit{Remark:} Any of the paths $\pi_i$ above can be viewed, according to 
Stone \cite{Stone1976} as a \textit{variation} of the geodesic $\overline{uv}$. Thus, 
in the setting of metric measure spaces \cite{Villani2009}, and for $PL$ manifolds 
\cite{Stone1976}, Eq. \ref{eq:haantjes-generalized} connects in the network context 
\textit{Jacobi fields} to Ricci curvature as in the classical case, see e.g. 
\cite{Jost2017}. 
\vspace{0.25cm}

We conclude this section by noting that both the version of the curvature for directed 
networks and weighted networks, can be extended, \textit{mutatis mutandis}, to this 
generalized definition.

\vspace{0.25cm}
\noindent \textit{Remark:} For the case of networks and simplicial complexes, it is 
useful to note that Eq. \ref{eq:haantjes-path} for a triangle $T = T(uvw)$ reduces 
to
\begin{equation}
\label{eq:haantjes-7}
\kappa^2_H(T) = \frac{d(u,v) + d(v,w) - d(u,w)}{(d(u,w))^3}\,.
\end{equation}
Thus, Haantjes curvature of triangles is closely related to two other measures, namely 
the \textit{excess} ${\rm exc}(T)$ and \textit{aspect ratio} ${\rm ar}(T)$, that are 
defined as follows:
\begin{equation}
\label{eq:excess}
{\rm exc}(T) = \max_{v \in \{u,v,w\}}{(d(u,v) + d(v,w) - d(u,w))}\,,
\end{equation}
and 
\begin{equation}
{\rm ar}(T) = \frac{exc(T)}{d(T)}\,,
\end{equation}
where $d(T)$ denotes the diameter of a triangle $T = T(uvw)$. There are strong 
connections between the excess, aspect ratio, and curvature. In particular, for the 
normalized Haantjes curvature introduced above, we have the following relation between 
the three notions:
\begin{equation}
\kappa_{H}^2(T(v)) = \frac{e(T(v))}{d^3}\,,
\end{equation}
that is
\begin{equation}
\kappa_{H}(T(v)) = \frac{\sqrt{{\rm ar}(T(v))}}{d}\,.
\end{equation}
Since the factor $\frac{1}{d}$ has the role of ensuring that, in the limit, the 
curvature of a triangle will have the dimensionality of the curvature at a point 
of a planar curve, the aspect ratio can be viewed as a (skewed), unnormalized 
version of the curvature, and Haantjes curvature can be viewed as a \textit{scaled} 
version of excess. Thus, curvature can be replaced by these surrogates, as the notion 
of scale is not of true import in many aspects of network characterization. Also, 
for the global understanding of the shape of networks, it is useful to compute, as 
is common in the manifold context, the \textit{maximal} excess and \textit{minimal} 
aspect ratio over all triangles in the network.
\vspace{0.25cm}

\vspace{0.25cm}
\noindent {\it Remark:} Note that Eqs. \ref{eq:haantjes-7} and \ref{eq:excess}
above again contain Gromov products \cite{Gromov2007}, that already appeared in 
the definition of Menger curvature.
\vspace{0.25cm}

\begin{figure*}
\centering
\includegraphics[width=8cm]{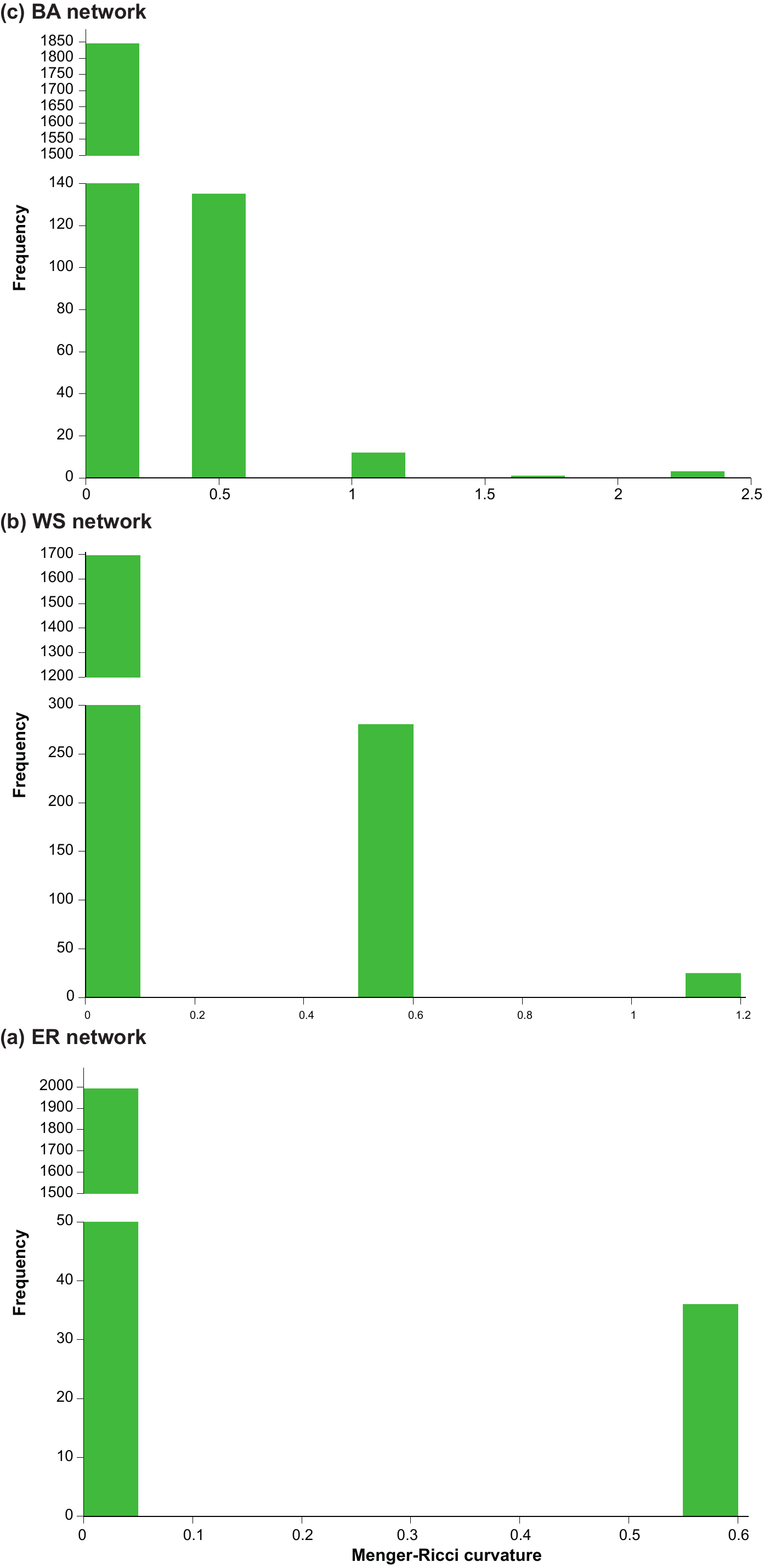}
\caption{Distribution of Menger-Ricci curvature in model networks. {\bf (a)}
ER network, {\bf (b)} WS network, and {\bf (c)} BA network with 1000 vertices 
and average degree 4.}
\end{figure*}

\begin{figure*}
\centering
\includegraphics[width=8cm]{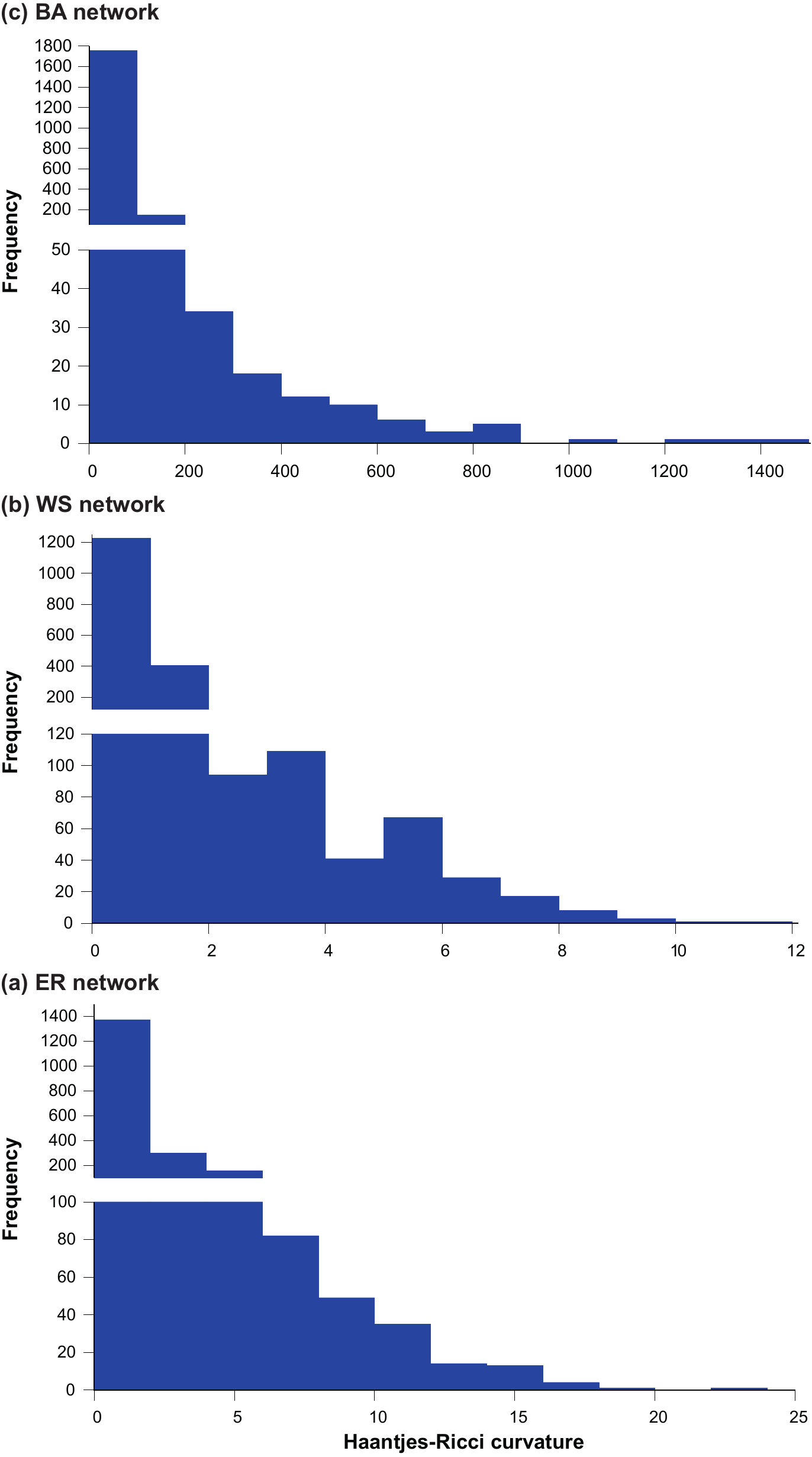}
\caption{Distribution of Haantjes-Ricci curvature in model networks. {\bf (a)}
ER network, {\bf (b)} WS network, and {\bf (c)} BA network with 1000 vertices 
and average degree 4. Here Haantjes-Ricci curvature is given by 
Eq. \ref{eq:haantjes-compute}.}
\end{figure*}

\section{Embedding Properties}

As noted in the Introduction, while studying complex networks, they are often 
considered for convenience to be equipped with the combinatorial metric, or 
they are considered to be topologically embedded in some model space and endowed 
with the induced (\textit{extrinsic}) metric, e.g. the hyperbolic one. However, 
at least in many instances, a more realistic approach would be to consider the 
networks to be endowed with an \textit{intrinsic} metric, obtained via some 
expressive (i.e., essential properties preserving) manner, from the weights
assigned to vertices and edges in the network. In this case, a natural question 
to ask is whether a (topological) embedding of the network into an ambient space 
will preserve the metric, or in other words, whether such an embedding is an 
\textit{isometric} one. 

Evidently, Haantjes-Ricci curvature, due to its purely metric definition, is 
\textit{geometric embedding invariant}: An isometric embedding will also preserve 
the curvature as is the case for Menger-Ricci curvature. Moreover, in the case 
where face weights are combinatorial, the sectional curvatures of faces, and hence, 
the Haantjes-Ricci curvatures are dependent solely on the metric structure of the 
edges, that is, ${\rm Ric}_H$ is a \textit{strong geometric embedding invariant} 
therefore face areas or weights are also preserved by an isometric embedding. However, 
in the case where face weights are general, the Haantjes-Ricci curvature is 
dependent on the face weights, as a consequence of Eq. \ref{eq:K-fitting} above, 
and we have the following result for general weighted networks.

\vspace{0.25cm}
\noindent \textbf{Proposition 4.1.} Let $(\mathcal{N},\mathcal{W})$ be a weighted 
network with general weights. Then, the Haantjes-Ricci curvature ${\rm Ric}_H$ is 
a geometric embedding invariant, however, not a strong geometric embedding one.
\vspace{0.25cm}

\noindent {\it Remark:} 
Any topological embedding that preserves edge and 2-face cell weights is a strong 
geometric embedding. In particular, for combinatorial (unweighted) networks, the 
Haantjes-Ricci curvature ${\rm Ric}_H$ is a strong geometric embedding invariant.
\vspace{0.25cm}


\begin{figure*}
\centering
\includegraphics[width=8cm]{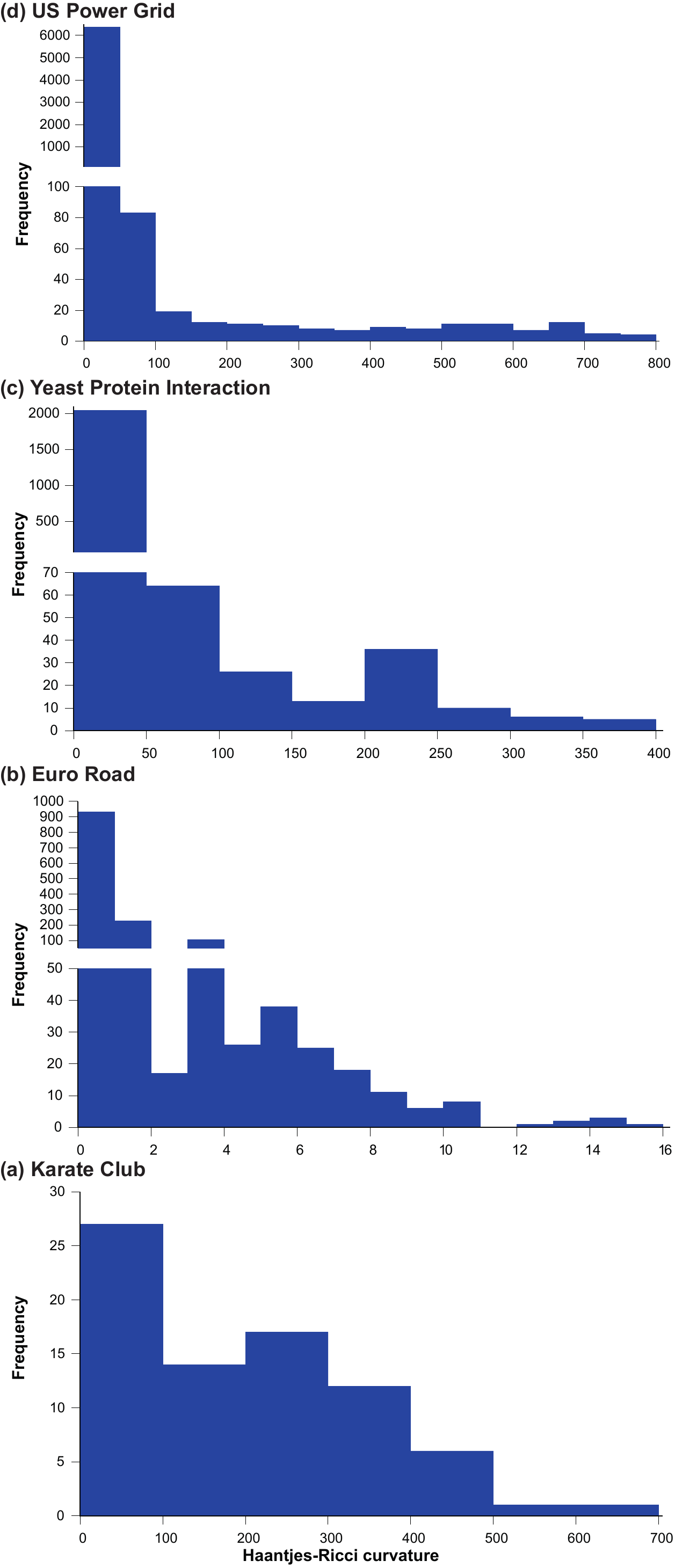}
\caption{Distribution of Haantjes-Ricci curvature in undirected real networks. 
{\bf (a)} Karate club, {\bf (b)} Euro road, {\bf (c)} Yeast protein interaction, 
and  {\bf (d)} US Power Grid.}
\end{figure*}

\begin{figure*}
\centering
\includegraphics[width=8cm]{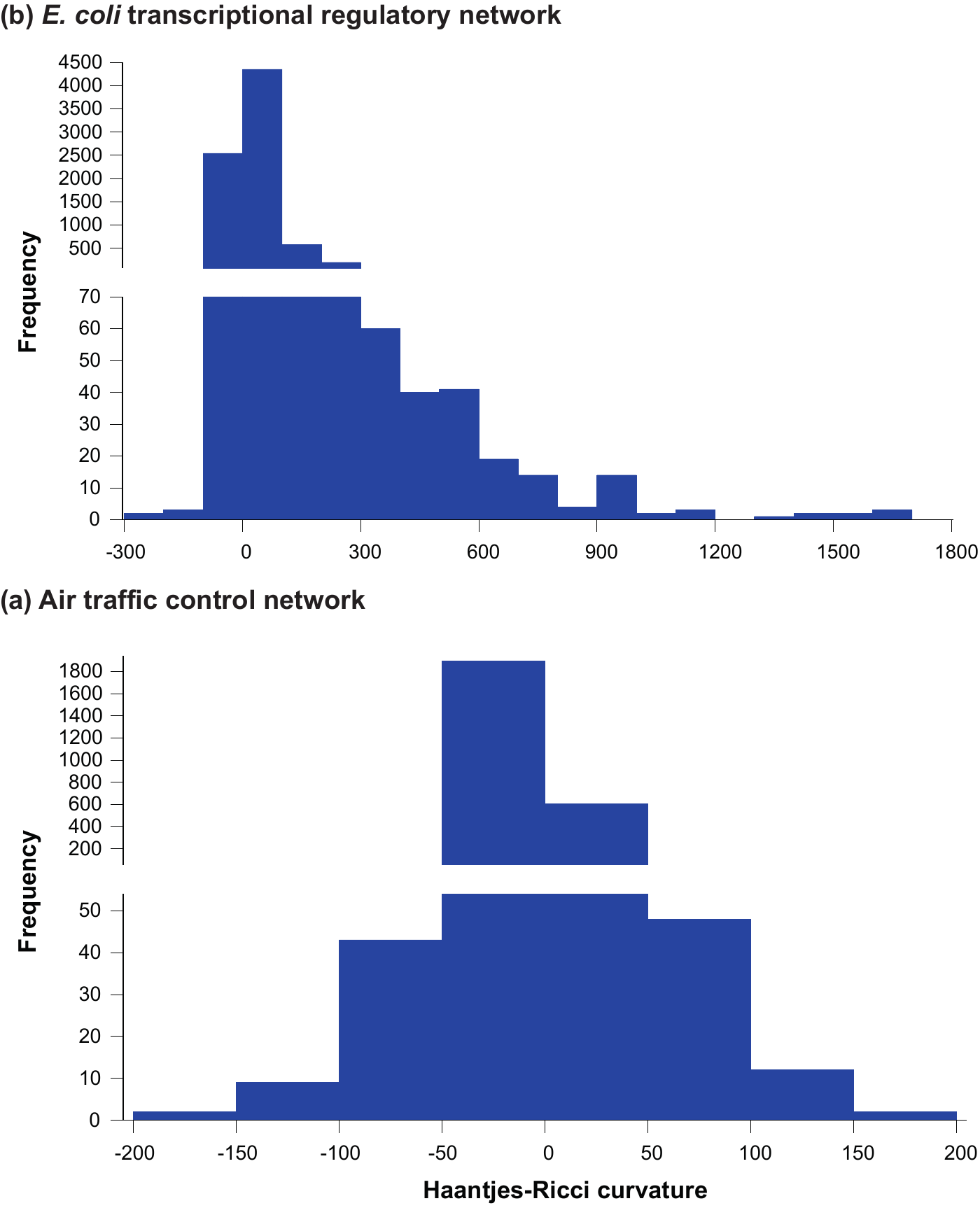}
\caption{Distribution of Haantjes-Ricci curvature in directed real networks. 
{\bf (a)} Air Traffic Control, and {\bf (b)} \textit{E. coli} transcriptional 
regulatory network. Note that, in this case of directed networks, the Haantjes-Ricci 
curvature of edges can be negative.}
\end{figure*}

\section{Application to model and real-world networks}

In this section, we explore the two notions of metric curvature introduced 
herein, in different types of model networks and real-world complex networks. 
Furthermore, we compare the two notions of metric curvature with existing 
network measures (including discrete notions of Ricci curvature) in both model 
and real-world complex networks. 

\subsection{Network datasets}

For this empirical exploration, we considered 3 models of undirected and 
unweighted complex networks, 4 undirected and unweighted real networks, and 2 
directed and unweighted real networks. 

The 3 models of undirected and unweighted complex networks considered here 
include the Erd\"{o}s-R{\'e}nyi (ER) model \cite{Erdos1961}, the Watts-Strogatz 
(WS) model \cite{Watts1998} and the Barab{\'a}si-Albert (BA) model \cite{Barabasi1999}. 
The ER model generates an ensemble of random graphs $G(n,p)$ where $n$ is the 
number of vertices and $p$ is the probability that any pair of vertices are 
connected via an edge in the network. The WS model generates small-world graphs 
with high clustering and small average path length. The parameters of the WS 
model are the number $n$ of vertices, the number $k$ of nearest neighbours to 
which each vertex is connected in the initial network, and the rewiring 
probability $\beta$ of edges in the initial network. The BA model generates 
scale-free graphs with power-law degree distribution. The parameters of the BA 
model are the number $n$ of vertices in the final network, the number $m_0$ of 
vertices in the initial network, and the number $m$ of existing vertices to which 
a new vertex is connected to at each step of this growing network model. Notably, 
the BA model uses a preferential attachment scheme whereby high-degree vertices 
have a higher chance of acquiring new edges than low-degree vertices at each step 
of this growing network model. 

The 4 undirected and unweighted real-world networks considered here are as 
follows. The Karate club network consists of 34 vertices representing members 
of the club, and 78 edges representing ties between pairs of members of the 
club, and this dataset was collected by Zachary in 1977 \cite{Zachary1977}. 
The Euro road network \cite{Subelj2011} consists of 1174 vertices representing 
cities in Europe, and 1417 edges representing roads in the international E-road 
network linking different cities. The Yeast protein interaction network 
\cite{Jeong2001} consists of 1870 vertices representing proteins in 
\textit{Saccharomyces cerevisiae}, and 2277 edges representing interactions 
between pairs of proteins. The US Power Grid network \cite{Leskovec2007} consists 
of 4941 vertices representing generators or transformers or substations in USA, 
and 6594 edges representing power supply lines linking them. The 2 directed 
and unweighted real-world networks considered here are as follows. The Air traffic 
control network \cite{Kunegis2013} consists of 1226 vertices representing airports 
and 2613 directed edges representing preferred routes between airports. The 
\textit{E. coli} transcriptional regulatory network (TRN) \cite{Salgado2013} 
consists of 3073 vertices representing genes and 7853 directed edges representing 
control of target gene expression through transcription factors.


\subsection{Distribution of metric curvatures in networks}

We have investigated the distribution of the two notions of metric curvature, 
namely, Menger-Ricci curvature and Haantjes-Ricci curvature in model networks 
considered here. Note that the Menger-Ricci curvature of edges in considered 
networks was computed using Eq. \ref{eq:Menger-Ricci} with Euclidean background 
geometry. In Figure 8, we display the distribution of Menger-Ricci curvature of 
edges in ER, WS and BA networks with 1000 vertices and average degree 4. From 
the figure, it is seen that the distribution of Menger-Ricci curvature of edges 
is broader in BA networks in comparison to WS networks which in turn is broader 
in comparison to ER networks, and this trend is preserved even when one compares 
ER, WS and BA networks with 1000 vertices and average degree 6 or 8 (data not 
shown).   

Note that the Haantjes-Ricci curvature of edges in considered networks was 
computed using Eq. \ref{eq:haantjes-compute}. Due to computational constraints, 
we only consider simple paths $\pi_i$ of length $\le 5$ between the two vertices 
anchoring any edge while computing its Haantjes-Ricci curvature using Eq. 
\ref{eq:haantjes-compute}. In Figure 9, we display the distribution of Haantjes-Ricci 
curvature of edges in ER, WS and BA networks with 1000 vertices and average degree 
4. From the figure, it is seen that the distribution of Haantjes-Ricci curvature 
of edges is very broad in BA networks in comparison to WS or ER networks, and 
this trend is preserved even when one compares ER, WS and BA networks with 1000 
vertices and average degree 6 or 8 (data not shown). 

In Figure 10, we display the distribution of Haantjes-Ricci curvature of edges in 
4 undirected and unweighted real networks considered here. In Figure 11, we display 
the distribution of Haantjes-Ricci curvature of directed edges in 2 directed and 
unweighted real networks considered here. As noted earlier, the Haantjes-Ricci 
curvature of directed edges can have negative value (unlike in the undirected and 
unweighted case).               


\subsection{Correlation between the two metric curvatures}

We next investigated the correlation between the two notions of metric curvature, 
namely, Menger-Ricci curvature and Haantjes-Ricci curvature in model and real 
networks considered here. Recall that, due to computational constraints, the 
computation of Haantjes-Ricci curvature of edges in considered networks using 
Eq. \ref{eq:haantjes-compute} was restricted to simple paths $\pi_i$ of maximum 
length $\le 5$ between the two vertices anchoring the edge under consideration.
Furthermore, we have also studied the correlation between Menger-Ricci curvature 
and Haantjes-Ricci curvature in model and real networks as a function of the 
maximum length of simple paths $\pi_i$ that are included in Eq. \ref{eq:haantjes-compute} 
while computing the Haantjes-Ricci curvature. 

In Table 2, we report the correlation between Menger-Ricci and Haantjes-Ricci 
curvature in model and real networks analyzed here. As expected, the Menger-Ricci 
curvature is perfectly correlated with Haantjes-Ricci curvature when the 
computation is restricted to paths of maximum length up to 2 or triangles 
(Table 2). However, the positive correlation between Menger-Ricci and 
Haantjes-Ricci curvature decreases with increase in the maximum length of paths 
accounted in the computation (Table 2). Specifically, Menger-Ricci and 
Haantjes-Ricci curvatures have minimal correlation in ER networks and moderate 
correlation in WS networks when paths of maximum length up to 5 are accounted 
in the computation. In contrast, there is high positive correlation between 
Menger-Ricci and Haantjes-Ricci curvatures in BA networks even when paths of 
maximum length up to 5 are accounted in the computation (Table 2). In real 
networks, we find moderate to high positive correlation between Menger-Ricci and 
Haantjes-Ricci curvatures even when paths of maximum length up to 5 are accounted 
in the computation (Table 2). Thus, one is conduced to the 
conclusion that the classical ER, WS and BA models have very little geometric 
content. This observation is corroborated  by the results in \cite{Weber2017}, 
that show that the distribution of Forman-Ricci curvature in real-life networks, 
has little to no resemblance to the one in the models above.  


\subsection{Comparison of the two metric curvatures with other network measures}

We next investigated the correlation of the two notions of metric curvature, 
Menger-Ricci and Haantjes-Ricci, with other network measures in model and real 
networks considered here. Specifically, we have studied the correlation with 
two other notions of Ricci curvature, (augmented) Forman-Ricci curvature and 
Ollivier-Ricci curvature \cite{Samal2018,Saucan2019a}, in networks. 
Specifically, the correlation we computed in the present paper is the classical 
Pearson correlation coefficient, which takes values in the range $[-1,+1]$.
In Table 3, we report the correlation between augmented Forman-Ricci curvature and 
Menger-Ricci or Haantjes-Ricci curvature in model and real networks analyzed 
here. We find no consistent trend in the correlation between augmented 
Forman-Ricci and Menger-Ricci curvature, or augmented Forman-Ricci and 
Haantjes-Ricci curvature in model and real networks analyzed here (Table 3). 
It is seen that there is high negative correlation between augmented Forman-Ricci 
and Haantjes-Ricci curvature in ER and BA networks, while there is high positive 
correlation between augmented Forman-Ricci and Menger-Ricci curvature in WS 
networks (Table 3). In Table 3, we also report the correlation between 
Ollivier-Ricci curvature and Menger-Ricci or Haantjes-Ricci curvature in model 
and real networks analyzed here. We again find no consistent trend in the 
correlation between Ollivier-Ricci and Menger-Ricci curvature, or Ollivier-Ricci 
and Haantjes-Ricci curvature in model and real networks analyzed here (Table 3).
Importantly, in recent work \cite{Samal2018,Saucan2019a}, it was shown that 
augmented Forman-Ricci and Ollivier-Ricci curvature have high positive correlation 
in model and real networks. These observations suggest that the two notions of 
metric curvature introduced herein capture different aspects of the network 
organization in comparison to discrete Ricci curvatures previously proposed for 
geometrical characterization of networks. This fact comes to further 
accentuate the fundamental observation that  it is not possible to find 
``the best'' discrete curvature, but rather the best one suited for a certain 
task, in a specific type of network. For example, given that Ollivier's 
curvature has the capability of predicting congestion in communication 
networks \cite{Wang2014,Wang2016}, thus it is excellent for the study of 
such network. However, it may be less suited for social or biological 
networks, since it is ``blind'' to cycles of length $\geq 6$.

Finally, we have also studied the correlation between edge betweenness centrality 
\cite{Freeman1977,Girvan2002} and Menger-Ricci or Haantjes-Ricci curvature in model 
and real networks (Table 4). We find no consistent trend in the correlation between 
edge betweenness centrality and Menger-Ricci curvature, or edge betweenness 
centrality and Haantjes-Ricci curvature in model and real networks analyzed here 
(Table 4). It is however seen that Menger-Ricci and Haantjes-Ricci curvature have 
moderate to high positive correlation with edge betweenness centrality in BA 
networks (Table 4). We remark that Menger-Ricci curvature, Haantjes-Ricci curvature, 
(augmented) Forman-Ricci curvature, Ollivier-Ricci curvature and edge betweenness 
centrality are edge-centric measures for analysis of networks, and thus, we have 
compared here these measures in model and real networks.


\begin{table}
\label{table:correlation}
\caption{Correlation between Menger-Ricci (MR) curvature and Haantjes-Ricci (HR) 
curvature in model and real networks. Here HR curvature is given 
by Eq. \ref{eq:haantjes-compute}.}
\begin{tabular}{l@{\hskip 0.43in}c@{\hskip 0.34in}c@{\hskip 0.34in}c@{\hskip 0.34in}c}
\hline \hline
{\bf Network}	& \multicolumn{4}{c}{{\bf MR versus HR curvature computed using}}	\\
				& \multicolumn{4}{c}{{\bf paths of maximum length up to}}			\\
				& {\bf 2}  		& {\bf 3} 		& {\bf 4} 	& {\bf 5} 				\\
\hline
\multicolumn{5}{l}{{\bf Undirected model networks}}\\
ER with average degree 4  	&	1	& 	0.31	&	0.14	&	0.08 	\\
ER with average degree 6	& 	1	& 	0.27	&	0.14	&	0.09 	\\
ER with average degree 8	& 	1	& 	0.24	&	0.12	&	0.11 	\\
WS with average degree 4 	&	1	&	0.75	&	0.55	&	0.38	\\
WS with average degree 6	&	1	&	0.77	&	0.63	&	0.44 	\\
WS with average degree 8	&	1	&	0.82	&	0.71	&	0.52	\\
BA with average degree 4 	&	1	&	0.59	&	0.52	&	0.51	\\
BA with average degree 6	&	1	&	0.71	&	0.71	&	0.71	\\
BA with average degree 8	&	1	&	0.88	&	0.87	&	0.86	\\
\hline
\multicolumn{5}{l}{{\bf Undirected real networks}}\\
Karate Club					&	1	& 	0.69	&	0.44	&	0.27 	\\
Euro Road					&	1	& 	0.51	&	0.37	&	0.32 	\\
Yeast Protein Interaction	&	1	&	0.83	&	0.74	&	0.64	\\
US Power Grid				&	1	&	0.82	&	0.67	&	0.58	\\
\hline
\multicolumn{5}{l}{{\bf Directed real networks}}\\
Air Traffic Control			&	1	& 	0.53	&	0.39	&	0.28	\\
{\it E. coli} TRN			&	1	&	0.83	&	0.75	&	0.71	\\
\hline \hline
\end{tabular}		
\end{table}	

\begin{table}
\label{table:ricci}
\caption{Comparison of augmented Forman-Ricci (AFR) curvature and Ollivier-Ricci 
(OR) curvature with Menger-Ricci (MR) curvature and Haantjes-Ricci (HR) curvature 
of edges in model and real networks. Here, HR curvature was computed by accounting 
for paths of maximum length up to 5 in Eq. \ref{eq:haantjes-compute}.}
\begin{tabular}{l@{\hskip 0.43in}c@{\hskip 0.25in}c@{\hskip 0.34in}c@{\hskip 0.25in}c}
\hline \hline
{\bf Network}	& \multicolumn{2}{c}{{\bf AFR curvature}} & \multicolumn{2}{c}{{\bf OR curvature}}\\
				& \multicolumn{2}{c}{{\bf versus}} & \multicolumn{2}{c}{{\bf versus}}\\
				& {\bf MR}  		& {\bf HR} 		& {\bf MR} 	& {\bf HR} \\
\hline
\multicolumn{5}{l}{{\bf Undirected model networks}}\\
ER with average degree 4  	&	0.04	& 	-0.55	&	0.08	&	-0.30 	\\
ER with average degree 6	& 	0.05	& 	-0.77	&	0.22	&	-0.31 	\\
ER with average degree 8	& 	0.06	& 	-0.86	&	0.29	&	-0.06 	\\
WS with average degree 4 	&	0.52	&	0.02	&	0.67	&	0.32	\\
WS with average degree 6	&	0.65	&	0.04	&	0.78	&	0.46 	\\
WS with average degree 8	&	0.70	&	0.06	&	0.80	&	0.52	\\
BA with average degree 4 	&	-0.37	&	-0.60	&	-0.10	&	-0.36	\\
BA with average degree 6	&	-0.44	&	-0.69	&	0.03	&	-0.13	\\
BA with average degree 8	&	-0.56	&	-0.69	&	0.17	&	0.09	\\
\hline
\multicolumn{5}{l}{{\bf Undirected real networks}}\\
Karate Club					&	0.45	& 	-0.22	&	0.13	&	-0.42 	\\
Euro Road					&	0.02	& 	-0.52	&	0.02	&	-0.25 	\\
Yeast Protein Interaction	&	0.12	&	-0.04	&	0.07	&	-0.10	\\
US Power Grid				&	0.11	&	-0.07	&	0.16	&	0.03	\\
\hline
\multicolumn{5}{l}{{\bf Directed real networks}}\\
Air Traffic Control			&	0.10	& 	0.13	&	-0.13	&	-0.02	\\
{\it E. coli} TRN			&	0.08	&	0.04	&	-0.06	&	-0.07	\\
\hline \hline
\end{tabular}		
\end{table}	

\begin{table}
\label{table:ebc}
\caption{Comparison of edge betweenness centrality with Menger-Ricci (MR) 
curvature and Haantjes-Ricci (HR) curvature of edges in model and real networks. 
Here, HR curvature was computed by accounting for paths of maximum length up 
to 5 in Eq. \ref{eq:haantjes-compute}.}
\begin{tabular}{l@{\hskip 0.43in}c@{\hskip 0.43in}c}
\hline \hline
{\bf Network}	& \multicolumn{2}{c}{{\bf Edge betweenness}} \\
				& \multicolumn{2}{c}{{\bf versus}} \\
				& {\bf MR}  		& {\bf HR} \\
\hline
\multicolumn{3}{l}{{\bf Undirected model networks}} \\
ER with average degree 4  	&	-0.03	& 	0.35	\\
ER with average degree 6	& 	-0.05	& 	0.51	\\
ER with average degree 8	& 	-0.06	& 	0.56	\\
WS with average degree 4 	&	-0.40	&	-0.20	\\
WS with average degree 6	&	-0.56	&	-0.24	\\
WS with average degree 8	&	-0.64	&	-0.27	\\
BA with average degree 4 	&	0.38	&	0.83	\\
BA with average degree 6	&	0.52	&	0.81	\\
BA with average degree 8	&	0.62	&	0.76	\\
\hline
\multicolumn{3}{l}{{\bf Undirected real networks}}\\
Karate Club					&	-0.20	& 	0.27	\\
Euro Road					&	-0.02	& 	0.04	\\
Yeast Protein Interaction	&	-0.07	&	0.05	\\
US Power Grid				&	-0.03	&	-0.02	\\
\hline
\multicolumn{3}{l}{{\bf Directed real networks}}\\
Air Traffic Control			&	-0.06	& 	-0.13	\\
{\it E. coli} TRN			&	-0.08	&	-0.06	\\
\hline \hline
\end{tabular}		
\end{table}	


\section{Conclusions and Future Outlook}

In previous work, notions of curvature for networks have been proposed, notably, 
Ollivier-Ricci curvature \cite{Ollivier2009} and Forman-Ricci curvature 
\cite{Forman2003}. Ollivier-Ricci curvature is cumbersome to compute in large 
networks while Forman-Ricci curvature is much less intuitive. Therefore, in 
search of simpler and more intuitive notions of curvature for networks and their 
higher-dimensional generalizations including simplicial and clique complexes, we 
have adapted here classical metric curvatures of curves, namely, introduced by 
Menger \cite{Menger1930} and Haantjes \cite{Haantjes1947}, for analysis of complex 
networks. In particular, we are able to define expressive metric Ricci curvatures 
for networks, both weighted and unweighted. Moreover, the simple yet elegant 
definitions of the metric curvatures introduced here are both computationally 
efficient (especially, those based on Menger curvature) and extremely versatile 
(especially, those derived from Haantjes curvature) from the perspective of 
application to complex networks. In fact, for the analysis of combinatorial 
(unweighted) polyhedral or simplicial complexes, the definitions of Haantjes-Ricci
curvature is shown here to be more expressive than the (augmented) Forman-Ricci
curvature as the Haantjes-Ricci curvature accounts for general $n$-gones and not 
only triangles like in the case of (augmented) Forman-Ricci curvature. 

Previously proposed metric curvature \cite{Saucan2012b} based on the \textit{Wald 
metric curvature} enables easy derivation of convergence results as well as proofs 
of theoretical results, such as a $PL$ version of the classical Bonnet-Myers 
theorem \cite{Gu2013}. Unfortunately, the metric curvatures proposed in 
\cite{Saucan2012b} are computationally expensive, and thus, impractical for 
analysis of large networks. This is in sharp contrast to the simplicity and 
efficiency of the metric curvatures for networks proposed in the present article.

We have also investigated model and real-world networks using the two notions of 
metric curvature introduced herein. It is seen that the distribution of Menger-Ricci 
curvature and Haantjes-Ricci curvature of edges is broader in scale-free BA networks 
in comparison to random ER networks or small-world WS networks. We also find a 
positive correlation between Menger-Ricci and Haantjes-Ricci curvature of edges 
in model and real-world networks, however, this correlation decreases with the 
increase in the maximum length of paths accounted in Eq. \ref{eq:haantjes-compute} 
for the computation of Haantjes-Ricci curvature. Thereafter, we have compared the 
Menger-Ricci or Haantjes-Ricci curvature with augmented Forman-Ricci or 
Ollivier-Ricci curvature in model and real networks, and no consistent trend in 
the correlation between Menger-Ricci or Haantjes-Ricci curvature with augmented 
Forman-Ricci or Ollivier-Ricci curvature was found in analyzed networks. Lastly, 
we also find no consistent trend in the correlation between edge betweenness 
centrality and Menger-Ricci or Haantjes-Ricci curvature in analyzed networks. Our
empirical results suggest that the metric curvatures introduced herein capture 
different aspects of the network organization in contrast to previously proposed 
edge-centric measures, such as Forman-Ricci curvature, Ollivier-Ricci curvature and
edge betweenness centrality, in complex networks. For instance, for combinatorial 
networks of bounded vertex degree, and in particular for finite networks, large 
Haantjes-Ricci curvature implies the existence of (arbitrarily) long simple cycles. 
Thus, the Haantjes-Ricci curvature, even if defined as a local invariant, has the 
potential of shedding light on the large scale topological structure of the network.

From a network application perspective, a natural extension is to develop algorithms  
based on Ricci curvatures introduced herein to detect clusters, modules or communities 
in real-world networks, and thereafter, compare the results obtained with algorithms 
for community detection in networks based on other notions of network curvature \cite{
Saucan2005,Sia2019,Ni2019}. Another application that imposes itself 
in view of the geodesic curvature role that Haantjes (and Menger) curvature plays 
in a vast, highly relevant class of metric spaces (see Appendix), is to use metric 
curvature for routing and hole detection ends.
Furthermore, another interesting exploratrion would be to 
investigate the correlation between the metric notions of curvature introduced herein 
and hyperbolic embeddings of networks. In particular, it is worthwhile to investigate 
the extent to which the curvatures can predict the values obtained using the inferred 
hyperbolic distances among vertices in embeddings considered in previous work 
\cite{Boguna2010}, and also study the extent to which these values are in concordance 
with the curvatures measured on the given network. Among the theoretical problems 
naturally arising from this work, one would like to first and foremost prove such 
results as a fitting analogue of the fundamental global Gauss-Bonnet theorem (which 
has important consequences in the study of long time evolution of networks 
\cite{Weber2018}), as well as a fitting Bonnet-Myers theorem.
Another important direction of theoretical study, that also has clear 
implications in real-life, large scale applications, is to determine the full extent 
of the local-to-global properties of the new invariants introduced in the present 
manuscript.


\subsection*{Acknowledgement}
The authors would like to thank the anonymous reviewers of the conference 
proceeding \cite{Saucan2019b} where some of these results were first reported. 
JJ and ES were partly supported by the German-Israeli Foundation Grant 
I-1514-304.6/2019. AS would like to thank Max Planck Society, Germany, for the 
award of a Max Planck Partner Group in Mathematical Biology.

\appendix
\section{Metric curvatures and classical curvature of curves} 
\label{sec:theory}

Both Menger curvature \cite{Menger1930} and Haantjes curvature \cite{Haantjes1947} 
were devised as generalizations to metric spaces of the classical notion of curvature 
for smooth plane curves. Therefore, it is only natural to ask whether they represent 
good approximations of the classical invariant. Indeed, as expected, this turns out 
to be the case for both metric curvatures. To state these results more formally we 
need the following definition:

\vspace{0.25cm}
\noindent \textbf{Definition A.1} Let $(M,d)$ be a metric space and let $p \in M$ 
be an accumulation point. We say that $M$ has \textit{Menger curvature} $\kappa_M(p)$ 
at $p$ iff for any $\varepsilon > 0$, there exists $\delta > 0$, such that for 
any triple of points $p_1, p_2, p_3$, satisfying $d(p,p_i) < \delta,\, i=1,2,3$, 
the following inequality holds: $|K_M(p_1,p_2,p_3) - \kappa_M(p)| < \varepsilon$.
\vspace{0.25cm}

While the corresponding definition for Haantjes curvature was essentially introduced 
in Eq. \ref{eq:haantjes-1}, we again state the same more formally here:

\vspace{0.25cm}
\noindent \textbf{Definition A.2} Let $(M,d)$ be a metric space and let $c: I=[0,1] 
\stackrel{\sim}{\rightarrow} M$ be a homeomorphism, and let $p,q,r \in c(I),\; q,r 
\neq p$. Denote by $\widehat{qr}$ the arc of $c(I)$ between $q$ and $r$, and by $qr$ 
line segment from $q$ to $r$. We say that the curve $c = c(I)$ has \textit{Haantjes} 
(or \textit{Finsler-Haantjes}) \textit{curvature} $\kappa_H(p)$ at the point $p$ iff:
\begin{equation}      
\label{eq:haantjes-2}
\kappa_H^2(p) = 24\lim_{q,r \rightarrow p}\frac{l(\widehat{qr})-d(q,r)}
{\big(l(\widehat{qr})\big)^3}\,\,,
\end{equation}
where $l(\widehat{qr})$ denotes the length, in the intrinsic metric induced by $d$, 
of $\widehat{qr}$. Note that the above definition is not the precise one used in this 
work.
\vspace{0.25cm}

For points or arcs where Haantjes curvature exists, we have $\frac{l(\widehat{qr})}{d(q,r)} 
\rightarrow 1$ as $d(q,r) \rightarrow 0$ (see \cite{Haantjes1947}), and therefore, 
$\kappa_{H}$ can be defined by (see e.g. \cite{Kay1980})
\begin{equation}                         
\label{eq:haantjes-3}
\kappa_{H}^2(p) = 24\lim_{q,r \rightarrow p}\frac{l(\widehat{qr})-d(q,r)}
{\big(d(q,r))\big)^3}\,\,.
\end{equation}
In the setting of graphs or networks, we prefer the above alternative form of the definition 
of Haantjes curvature which is even more intuitive and simpler to use.

As expected, Haantjes curvature for smooth curves in the Euclidean plane (or space) coincides 
with the standard (differential) notion. More precisely, we have the following result, see 
\cite{Haantjes1947}.

\vspace{0.25cm}
\noindent \textbf{Theorem A.3} Let  $c \in \mathcal{C}^3$ be smooth curve in $\mathbb{R}^3$ 
and let $p \in c$ be a regular point. Then the metric curvature $\kappa_{H}(p)$ exists and equals 
the classical curvature of $\gamma$ at $p$.
\vspace{0.25cm}

A similar result also holds for Menger curvature. However, we will arrive at the result in an 
indirect fashion which will enable us to establish the connection between the two types of metric 
curvature used here. Firstly, note that although the Haantjes curvature exhibits additional 
flexibility in the network context in comparison to the Menger curvature, the formal definition of 
Haantjes curvature in the setting of metric spaces is a more restricted notion than the Menger 
curvature, since it can be employed only for rectifiable curves. However, the following theorem 
due to \cite{Pauc1936} holds:

\vspace{0.25cm}
\noindent \textbf{Theorem A.4} Let $(X,d)$ be a metric continuum, and consider $p \in X$. If 
$\kappa_M(p)$ exists, then $X$ is a rectifiable arc in a neighbourhood of $p$.
\vspace{0.25cm}

\vspace{0.25cm}
\noindent \textbf{Corollary A.5} Let $(X,d)$ be a metric arc. If $\kappa_M(p)$ exists at all 
points $p \in X$, then $X$ is rectifiable.
\vspace{0.25cm}

In consequence, the existence of one of the considered metric curvatures implies the existence 
of the other. In fact, the two definitions coincide whenever they both are applicable:

\vspace{0.25cm}
\noindent \textbf{Theorem A.6} \cite{Haantjes1947} Let $c$ be a rectifiable arc in a metric 
space $(M,d)$, and let $p \in c$. If $\kappa_M$  and $\kappa_{H}$ exist, then they are equal.
\vspace{0.25cm}

Furthermore, it turns out that another important property of classical curvature, 
namely characterization of geodesics via {\it geodesic curvature}  (see, e.g. 
\cite{Jost2017}) is also captured, at least for a large class of metric spaces, 
by Haantjes curvature (hence, due to Theorem A.6 above, by Menger's curvature 
as well). More precisely, we have the following result:

\vspace{0.25cm}
\noindent \textbf{Theorem A.7} \cite{Haantjes1947} 	Let $(X,d)$ be a metric 
space and let $\gamma \subset X$ be a {\em Ptolemaic} arc (in the induced metric). 
$\gamma$ is a geodesic segment if $\kappa_H(p) = 0$, for any $p \in \gamma$.
\vspace{0.25cm}

While we do not bring here the technical definition of Ptolemaic metric 
spaces, since this would expand the present text in an inordinate manner 
and we only refer the reader to \cite{Blumenthal1953,Blumenthal1970}. We 
wish to underline that Euclidean and Hyperbolic spaces are Ptolemaic, 
thus one can apply the theorem above to networks embedded in such spaces.





\end{document}